\documentclass[11pt]{amsart}

\usepackage[T1]{fontenc}
\usepackage[utf8]{inputenc}
\usepackage{fouriernc}    
\usepackage{microtype}

\usepackage[
  top    = 1.35in,
  bottom = 1.35in,
  inner  = 1.25in,
  outer  = 1.25in,
  footskip = 30pt
]{geometry}

\usepackage{titlesec}

\titleformat{\section}
  [block]
  {\normalfont\scshape\centering}
  {\S\,\Roman{section}.}
  {0.6em}
  {}
  [\vspace{3pt}{\titlerule[0.4pt]}\vspace{-2pt}]

\titlespacing*{\section}{0pt}{2.4ex plus .4ex}{1.6ex plus .2ex}

\titleformat{\subsection}
  [runin]
  {\normalfont\scshape}
  {\thesubsection.}
  {0.5em}
  {}
  [.\enspace]

\titlespacing*{\subsection}{0pt}{1.4ex plus .3ex}{0pt}

\renewcommand{\thesection}{\Roman{section}}
\renewcommand{\thesubsection}{\thesection.\arabic{subsection}}

\usepackage{amsmath, amssymb, amsfonts, amsthm, mathtools}
\numberwithin{equation}{section}
 \usepackage{booktabs}
 \usepackage{array}
 \usepackage{tabularx}
 \usepackage{makecell}

\newtheoremstyle{classical}
  {7pt}{5pt}{\itshape}{\parindent}{\scshape}{.}{ }{}

\newtheoremstyle{classicaldef}
  {7pt}{5pt}{\normalfont}{\parindent}{\scshape}{.}{ }{}

\theoremstyle{classical}
\newtheorem{theorem}{Theorem}[section]
\newtheorem{proposition}[theorem]{Proposition}
\newtheorem{lemma}[theorem]{Lemma}

\theoremstyle{classicaldef}
\newtheorem{definition}[theorem]{Definition}

\newtheorem{remark}[theorem]{Remark}

\usepackage{graphicx}
\usepackage{booktabs}
\usepackage{array}
\usepackage{subcaption}
\usepackage{tikz}
\usetikzlibrary{
  fit, backgrounds, arrows.meta,
  positioning, calc, decorations.pathreplacing
}
\usepackage{tikz-cd}
\usepackage{float}

\usepackage{enumitem}
\setlist{leftmargin=*, itemsep=2pt, topsep=4pt}

\usepackage[square,numbers,sort&compress]{natbib}
\usepackage[hidelinks]{hyperref}
\usepackage[nameinlink,noabbrev]{cleveref}

\begin{document}

\title{On $ \mathrm{U(1)}^{2} $-Invariant Cayley Submanifolds in $ \mathbb{R}^{8} $}
\author{Mia S.L.\ Beard \\ \textit{Mathematical Institute, University of Oxford}}

\date{\today}
\maketitle
 
\begin{abstract}

   We study $ \mathrm{U}(1)^{2} $-invariant Cayley submanifolds in $ \mathbb{R}^{8} $, under the action with weights $ (e^{i \theta_{1}}, e^{i \theta_{2}}, e^{-i(\theta_{1} + \theta_{2})}, 1) $ on coordinates $ (z_{1}, z_{2}, z_{3}, z_{4}) \in \mathbb{C}^{4} $. Locally, we can prove that $ N $ can be written as a graph of functions $ u, v, a_{1}, a_{2}: \mathbb{R}^{2} \rightarrow \mathbb{R} $ (where $ a_{3}= -a_{1}-a_{2} $), where unlike in the special Lagrangian setting, the moment-map parameters are not required to be constant, which gives rise to additional degrees of freedom in the reduced equations. We prove that $ N $ is Cayley if and only if a non-linear, Cauchy-Riemann system of partial differential equations in $ u, v, a_{1}, a_{2} $ holds, with ellipticity away from its degeneracy locus. 
    
    We investigate the special Lagrangian and coassociative reductions of the system, obtaining potential functions in both cases. For general $ \mathrm{U}(1)^{2} $-invariant Cayley $ 4 $-folds, we derive almost-conservation laws (with clear obstruction terms) involving natural $ S_{3} $-invariant and alternating combinations of the moment-map parameters. These identities exhibit additional source terms, distinguishing the general Cayley system from its special Lagrangian and coassociative reductions. 
    
\end{abstract}
 
\vspace{1.5em}

\section*{Introduction}

A calibration on a Riemannian manifold $ (X, g) $ is a closed differential $ k $-form $ \varphi $, whose restriction to every oriented tangent $ k $-plane is bounded above by the corresponding volume form. An oriented $ k $-dimensional submanifold on which equality holds is said to be calibrated. Compact calibrated submanifolds are volume-minimising in their homology classes, and form a distinguished class of minimal submanifolds arising naturally in geometries with special holonomy.

Several interesting families of calibrated submanifolds come from consideration of special holonomy groups.

Of particular interest are Cayley $ 4 $-folds, the calibrated submanifolds which are associated with $ \mathrm{Spin}(7) $-geometry. Identifying $ \mathbb{R}^{8} $ with $ \mathbb{C}^{4} $, and equipping the latter with its standard Euclidean metric, K\"{a}hler form $ \omega $, and holomorphic volume form $ \Omega = dz_{1} \wedge dz_{2} \wedge dz_{3} \wedge dz_{4} $, the standard Cayley calibration is

\begin{equation*}
    \Phi = \frac{1}{2}\omega \wedge \omega + \mathrm{Re}(\Omega).
\end{equation*}

The defining equations for Cayley submanifolds reveal a unifying relationship between several calibrated geometries. Complex surfaces and special Lagrangian $ 4 $-folds in $ \mathbb{C}^{4} $ are Cayley, while Cayley $ 4 $-folds contained in an affine hyperplane $ \mathbb{R}^{7} \subset \mathbb{R}^{8} $ are coassociative with respect to the Hodge dual of the induced $ G_{2} $-structure. Cayley geometries are therefore a common setting in which these different classes of calibrated submanifolds may be studied.

A natural approach to constructing calibrated submanifolds is to impose symmetry, reducing the calibration condition to a system of differential equations on a lower-dimensional space. For $ \mathrm{U(1)} $-invariant special Lagrangian submanifolds in $ \mathbb{C}^{3} $, this approach has been developed extensively by Joyce \cite{joyce2005u1, joyce2005u2, joyce2005u3}, and extended by this paper's author in \cite{MyPaper} to $ \mathrm{U(1)}^{n-2} $-invariant special Lagrangian $ n $-folds in $ \mathbb{C}^{n} $. A central feature of these constructions is that the moment map associated to the symmetry group is constant along each connected Lagrangian submanifold.

The Cayley condition is more general, in particular: the moment-map parameters are no longer required to be constant. This observation motivates the present paper, in which we develop a system of partial differential equations characterising a class of $ \mathrm{U(1)}^{2} $-invariant Cayley $ 4 $-folds in $ \mathbb{C}^{4} $, where there are two additional partial differential equations compared to the special Lagrangian construct, representing the fact that the moment-map parameters themselves now vary over the two-dimensional reduced domain.

Consider the $ \mathrm{U(1)}^{2} $-action in coordinates $ \mathbb{C}^{4} $

\begin{equation*}
    (e^{i \theta_{1}}, e^{i \theta_{2}}) \cdot (z_{1}, z_{2}, z_{3}, z_{4}) = \bigl(e^{i \theta_{1}}z_{1}, e^{i \theta_{2}}z_{2}, e^{-i(\theta_{1}+\theta_{2})}z_{3}, z_{4} \bigr).
\end{equation*}

This action preserves the standard Cayley calibration. The product $ iz_{1}z_{2}z_{3} $ and the coordinate $ z_{4} $ are invariant, while the differences $ |z_{1}|^{2} - |z_{3}|^{2} $ and $ |z_{2}|^{2} - |z_{3}|^{2} $ determine the moment-map parameters, up to a conventional normalisation.

Let $ \Sigma \subset \mathbb{R}^{2} $ be a domain with coordinates $ (x, y) $, and let $ u, v, a_{1}, a_{2}: \Sigma \rightarrow \mathbb{R} $ be smooth functions. Set $ a_{3} = -a_{1} - a_{2} $, and define a $ \mathrm{U(1)}^{2} $-invariant subset $ N \subset \mathbb{C}^{4} $ by

\begin{equation*}
\begin{aligned}
    N := \Big\{ z \in \mathbb{C}^{4} \;\Big|\;
    &|z_{1}|^{2} - |z_{3}|^{2} = a_{1}(x,y) - a_{3}(x,y), \\
    &|z_{2}|^{2} - |z_{3}|^{2} = a_{2}(x,y) - a_{3}(x,y), \\
    &iz_{1}z_{2}z_{3} = v(x,y) + iy, \\
    &z_{4} = x+iu(x,y), \quad (x,y) \in \Sigma \Big\}.
\end{aligned}
\end{equation*}

The three squared radii may be expressed as $ |z_{j}|^{2} = w+a_{j} $, where $ w = w(x, y) $ is determined by $ P(w) = (w+a_{1})(w+a_{2})(w+a_{3}) = v^{2}+y^{2} $, subject to the admissibility condition $ w \geq -\min\{a_{1}, a_{2}, a_{3} \} $. At each point, this condition determines a unique admissible root. 

The key result of this paper (presented in \S \ref{SECT:ConstructionOfCayleyFourFolds}) establishes that if we work on the domain $ \Sigma^{*} $, by which we mean $ \Sigma $ minus the set where the $ \mathrm{U(1)}^{2} $-orbit degenerates, with $ u, v, a_{1}, a_{2} $ as smooth on $ \Sigma^{*} $, then $ N $ is Cayley if and only if: 

\begin{equation*}
    \begin{cases}
        u_{x} = v_{y}, \\
        v_{x} = -P'(w)u_{y}, \\
        a_{1, x} = -((a_{1}+a_{2}+w)a_{1, y} + (a_{1}+2w)a_{2, y}), \\
        a_{2, x} = (a_{2}+2w)a_{1, y} + (a_{1}+a_{2}+w)a_{2, y}.
    \end{cases}
\end{equation*}

The reduced system is elliptic on the regular locus, where $ P'(w) > 0 $.

The $ \mathrm{U(1)}^{2} $ action fails to be free at $ z $ precisely when $ v(x, y) = y = 0 $ and the minimum of $ \allowbreak{\{a_{1}(x, y), a_{2}(x, y), a_{3}(x, y) \}} $ is attained with multiplicity at least two. 

The system naturally separates into three regimes. 

If $ u $ is constant, then the first two equations force $ v $ to be constant on each connected component of the regular domain. The resulting submanifold lies in the affine hyperplane $ H_{u_{0}}
    =\{z\in\mathbb{C}^{4}:\operatorname{Im}z_{4}=u_{0}\}
    \cong\mathbb{R}^{7} $, and is coassociative with respect to its reduced $ G_{2} $-structure.

These three regimes have different analytic structures. In the full Cayley case, working with the natural symmetric polynomials yields no obvious potential function, although we do observe an `almost conservation law', with an obstruction term measuring the failure to admit a conservation law. By contrast, the special Lagrangian equations admit a potential function, which underlies the Dirichlet existence and uniqueness theory developed in \cite{joyce_special_2002} and discussed in higher generality in \cite{MyPaper}. Perhaps more unexpectedly, the coassociative equations also admit a potential function, reducing the system to a single quasilinear second-order elliptic partial differential equation. This appears to warrant further study of the existence and uniqueness of $ \mathrm{U(1)}^{2} $-invariant coassociative submanifolds.

\newcolumntype{Y}[1]{>{\raggedright\arraybackslash}p{#1\linewidth}}
\begin{table}[htbp]
    \centering
    \small
    \renewcommand{\arraystretch}{1.35}
    \caption{The different regimes of our construction.}
    \label{tab:regimes}
    \begin{tabularx}{\textwidth}{@{} l l c Y{0.22} Y{0.20} @{}}
        \toprule
        \textbf{Sector} & \textbf{Conditions} & \textbf{Reduced} & \textbf{Potential} & \textbf{Calibration} \\
        & & \textbf{variables} & & \\
        \midrule
        Special Lagrangian &
        $a_1, a_2, a_3$ constant &
        $u, v$ &
        $f_x = v, \quad f_y = u$ &
        Special Lagrangian in $\mathbb{C}^4$ \\[4pt]
        Coassociative &
        $u = u_0,\ v = v_0$ constant &
        $a_1, a_2, a_3$ &
        \makecell[l]{$\begin{aligned}
    g_x &= -\tfrac{1}{3}(a_1-a_2)(a_2-a_3) \\
        &\quad \times(a_3-a_1) \\
    g_y &= a_1a_2 + a_2a_3 + a_3a_1 \\
        &\quad + \tfrac{3}{2}w^2
\end{aligned}$} &
        Coassociative in $H_{u_0} \cong \mathbb{R}^7$ \\[10pt]
        Genuine Cayley &
        Both vary &
        $u, v, a_1, a_2, a_3$ &
        Obstruction &
        Cayley in $\mathbb{C}^4$ \\
        \bottomrule
    \end{tabularx}
\end{table}

\section{Prerequisites in Calibrated Geometry}\label{SECT:Prerequisites}

\subsection{Rudiments of Calibrated Geometry}

We begin by reviewing the basic definitions and results in calibrated geometry which are applicable to this paper, using the standard references of \cite{harvey_calibrated_1982, joyce2000compact, joyce2007riemannian}.

\begin{definition}[Oriented tangent $ k $-planes and volume forms]

Let $ (M, g) $ be a Riemannian manifold. An \textit{oriented tangent $ k $-plane} at a point $ p \in M $ is a $ k $-dimensional vector subspace $ V \subset T_{p}M $, together with a choice of orientation. The restriction of $ g $ to $ V $ is an inner product, and hence determines a Euclidean volume form $ \mathrm{Vol}_{V} \in \Lambda^{k} V^{*} $. Equivalently, if $ (e^{1}, \dots, e^{k}) $, then

\begin{equation*}
    \mathrm{Vol}_{V} = e^{1} \wedge \cdots \wedge e^{k},
\end{equation*}

which spans the one-dimensional vector space $ \Lambda^{k}V^{*} $.
    
\end{definition}

\begin{definition}[Calibrations and calibrated submanifolds]

Let $ \varphi \in \Omega^{k}(M) $ be a closed $ k $-form. We say that $ \varphi $ is a \textit{calibration} on $ (M, g) $ if, for every oriented tangent $ k $-plane $ V \subset T_{p}M $,

\begin{equation}\label{EQ:Calibration}
    \varphi|_{V} \leq \mathrm{Vol}_{V}.
\end{equation}

Equivalently, for every positively oriented orthonormal basis $ (e_{1}, \dots, e_{k}) $ of $ V $,

\begin{equation*}
    \varphi_{p}(e_{1}, \dots, e_{k}) \leq 1.
\end{equation*}

Let $ N \subset M $ be an oriented $ k $-dimensional submanifold. We say that $ N $ is \textit{calibrated by $ \varphi $} if $ \varphi|_{T_{p}N} = \mathrm{Vol}|_{T_{p}N} $ for every $ p \in N $. In other words, the tangent planes of $ N $ are precisely those on which $ \varphi $ attains its maximal positive value.
    
\end{definition}

The fundamental theorem of calibrated geometry is that all calibrated submanifolds minimise volume \cite{harvey_calibrated_1982}.

\begin{theorem}

    All compact calibrated submanifolds are volume minimising in their homology class.
   
\end{theorem}

Many examples of calibrated manifolds are generated by considering holonomy groups. Let $ G \subset O(n) $ be a possible holonomy group of a Riemann metric, where we could, for example, take $ G $ to be $ U(m), SU(m), Sp(m), Sp(m)Sp(1), G_{2} $, or $ \mathrm{Spin}(7) $. In this paper, we shall consider the following.

\begin{itemize}

    \item Let $ G = SU(m) \subset O(2m) $. A manifold $ (M, g) $ with holonomy $ SU(m) $ is called a \textit{Calabi-Yau manifold}. This comes equipped with the K\"{a}hler form $ \omega $ and the holomorphic volume form $ \theta $ The calibrated submanifolds with this reduced holonomy are called \textit{special Lagrangian submanifolds}.
 
    \item The group $ G = G_{2} \subset O(7) $ preserves a $ 3 $-form $ \varphi_{0} $ and a $ 4 $-form $ * \varphi_{0} $ on $ \mathbb{R}^{7} $. A manifold $ (M, g) $ with holonomy $ G_{2} $ has two different forms: $ \varphi $, which is an associative $ 3 $-form, and $ * \varphi $, which is a coassociative $ 4 $-form. The calibrated submanifolds with reduced holonomy $ \varphi $ are called associative $ 3 $-folds, and those with reduced holonomy $ * \varphi $ are called coassociative $ 4 $-folds.

    \item The group $ G = \mathrm{Spin}(7) \subset O(8) $ preserves a $ 4 $-form $ \Omega_{0} $ on $ \mathbb{R}^{8} $, thus $ (M^{8}, g) $  with holonomy $ \mathrm{Spin}(7) $ has a $ 4 $-form $ \Omega $, which is a calibration. We call submanifolds which are calibrated with respect to $ \Omega $ \textit{Cayley $ 4 $-folds}.
    
\end{itemize}

\subsection{\texorpdfstring{Special Lagrangian Submanifolds in $\mathbb{C}^{4}$}{Special Lagrangian Submanifolds in C4}}

\begin{definition}[Special Lagrangian submanifolds in $ \mathbb{C}^{4} $]

Let $ (z_{1}, z_{2}, z_{3}, z_{4}) \in \mathbb{C}^{4}$, which is equipped with the Riemannian metric $ g = |dz_{1}|^{2} + |dz_{2}|^{2} + |dz_{3}|^{2} + |dz_{4}|^{2} $, the K\"{a}hler form $ \omega = \tfrac{i}{2}(dz_{1} \wedge d \bar{z}_{1} + dz_{2} \wedge d \bar{z}_{2} + dz_{3} \wedge d \bar{z}_{3} + dz_{4} \wedge d \bar{z}_{4}) $, and the holomorphic volume form $ \Omega = dz_{1} \wedge dz_{2} \wedge dz_{3} \wedge dz_{4} $. 

Let $ L \subset \mathbb{C}^{4} $ be an oriented real four-dimensional submanifold. We say that $ L $ is a \textit{special Lagrangian submanifold} of phase zero if it is calibrated with respect to $ \mathrm{Re}(\Omega) $, that is,

\begin{equation*}
    \mathrm{Re}(\Omega)|_{L} = \mathrm{Vol}_{L}.
\end{equation*}

Equivalently, a real four-dimensional submanifold $ L \subset \mathbb{C}^{4} $ admits an orientation making it special Lagrangian of phase zero if and only if

\begin{equation*}
    \omega|_{L} \equiv 0, \qquad \mathrm{Im}(\Omega)|_{L} \equiv 0,
\end{equation*}

with the orientation chosen so that $ \mathrm{Re}(\Omega)|_{L} = \mathrm{Vol}_{L} $.
    
\end{definition}

There is well-established theory surrounding special Lagrangian $ n $-folds. Of particular relevance here is Joyce's paper \cite{joyce2005u1}, which establishes a construction for $ \mathrm{U(1)} $-invariant special Lagrangian $ 3 $-folds in $ \mathbb{C}^{3} $, whereby a $ \mathrm{U(1)} $-invariant submanifold $ N $, which can be written as a graph of functions $ u, v: \mathbb{R}^{2} \to \mathbb{R} $, is realised as a special Lagrangian $ 3 $-fold if and only if it satisfies a system of non-linear Cauchy-Riemann equations which depend on $ a $, where $ a $ is defined by $ |z_{1}|^{2} - |z_{2}|^{2} = 2a $. If $ a \neq 0 $, then $ u, v $ are smooth and $ N $ is a non-singular special Lagrangian $ 3 $-fold, but if $ a = 0 $, then there may be points of the form $ (x, y) $, where $ u, v $ are not differentiable, which corresponds to singular points of $ N $. In the non-singular case, the existence and uniqueness of a large class of non-singular $ \mathrm{U(1)} $-invariant special Lagrangian $ 3 $-folds follows intuitively. In Joyce's further papers \cite{joyce2005u2, joyce2005u3}, it is recognised that in the case of $ a = 0 $, the Cauchy-Riemann equations are not always elliptic, leading to technical analytic issues. Proving existence and uniqueness in this case is proven using a priori estimates for derivatives of solutions, using these to show existence and uniqueness of weak solutions $ u, v $ to two Dirichlet problems when $ a = 0 $, which are continuous and weakly differentiable. In studying the singularities, it is found that, under mild conditions, the singularities are isolated, and have a multiplicity $ n > 0 $, and one of two types. Further, the author of this paper formulated \cite{MyPaper}, which generalises this result to the construction of $ \mathrm{U(1)}^{n-2} $-invariant special Lagrangian $ n $-folds in $ \mathbb{C}^{n} $. There is thus well-established theory surrounding the existence and uniqueness of these objects, as well as singularity analysis.

\subsection{\texorpdfstring{Coassociative Submanifolds in $ \mathbb{R}^{7} $}{Coassociative Submanifolds in R⁷}}

There are two types of calibrated submanifolds in $ \mathbb{R}^{7} $: associative $ 3 $-folds and coassociative $ 4 $-folds. 

\begin{definition}[Associative $ 3 $-Folds]

Take coordinates $ (x_{1}, \dots, x_{7}) \in \mathbb{R}^{7} $, a metric $ g_{0} = \rm{d}x_{1}^{2} + \cdots + \rm{d}x_{7}^{2} $, and a $ 3 $-form $ \varphi_{0} $ on $ \mathbb{R}^{7} $ by

\begin{equation}\label{EQ:AssociativeSubmanifoldsCalibration}
    \varphi_{0} = \rm{d}x_{123} + \rm{d}x_{145} + \rm{d}x_{167} + \rm{d}x_{246} - \rm{d}x_{257} - \rm{d}x_{347} - \rm{d}x_{356},
\end{equation}

where $ \rm{d}x_{ijk} = \rm{d}x_{i} \wedge \rm{d}x_{j} \wedge \rm{d}x_{k} $.

An \textit{associative $ 3 $-fold} in $ \mathbb{R}^{7} $ is an oriented three-dimensional submanifold of $ \mathbb{R}^{7} $ which is calibrated with respect to $ \varphi_{0} $.
    
\end{definition}

\begin{definition}[Coassociative $ 4 $-Folds]

Take coordinates $ (x_{1}, \cdots, x_{7}) \in \mathbb{R}^{7} $, a metric $ g_{0} = dx_{1}^{2} + \cdots + dx_{7}^{2} $, and define a $ 4 $-form $ * \varphi_{0} $ on $ \mathbb{R}^{7} $, where

\begin{equation}\label{EQ:CoassociativeSubmanifoldsCalibration}
    * \varphi_{0} = \rm{d}x_{4567} + \rm{d}x_{2367} + \rm{d}x_{2345} + \rm{d}x_{1357} - \rm{d}x_{1346} - \rm{d}x_{1256} - \rm{d}x_{1247}.
\end{equation}

A \textit{coassociative $ 4 $-fold} in $ \mathbb{R}^{7} $ is an oriented four-dimensional submanifold of $ \mathbb{R}^{7} $ calibrated with respect to $ * \varphi_{0} $.
    
\end{definition}

The primary research which has been carried out into the construction of coassociative submanifolds is Lotay's paper \cite{lotay2007calibrated},which establishes a technique for constructive coassociative submanifolds of Euclidean space with symmetry, using first-order ODEs (there are, however, no obvious solutions to these ODEs).

\subsection{Cayley Submanifolds in $ \mathbb{R}^{8} $}

\begin{definition}[Cayley $ 4 $-folds in $ \mathbb{R}^{8} $]

Take coordinates $ (x_{1}, \dots, x_{8}) \in \mathbb{R}^{8} $, and define the metric $ g_{0} = dx_{1}^{2} + \cdots + dx^{8} $, and a $ 4 $-form $ \Omega $ by

\begin{equation}
\begin{aligned}
    \Omega = \rm{d}x_{1234}+\rm{d}x_{1256} + \rm{d}x_{1278} + \rm{d}x_{1357} - \rm{d}x_{1368} - \rm{d}x_{1458} - \rm{d}x_{1467} \\
    - \rm{d}x_{2358} - \rm{d}x_{2367} - \rm{d}x_{2457} + \rm{d}x_{2468} + \rm{d}x_{3456} + \rm{d}x_{3478} + \rm{d}x_{5678}.
\end{aligned}
\end{equation}

A \textit{Cayley $ 4 $-fold} in $ \mathbb{R}^{8} $ is an oriented $ 4 $-dimensional submanifold of $ \mathbb{R}^{8} $ which is calibrated with respect to $ \Omega $.
    
\end{definition}

Several papers construct explicit examples of Cayley $ 4 $-folds in $ \mathbb{R}^{8} $: Gu and Pries \cite{ExamplesOfCayley4Manifolds}, Ionel, Karigiannis, and Min-oo \cite{BundleConstructionInR7AndR8}, and Lotay, in two seperate papers \cite{ lotay20062, lotay2007calibrated}. The geometry of Cayley submanifolds of Calabi-Yau four-folds has also been studied, by Ghigi \cite{ghigi2000generalization}, and by Salavessa and do Vale \cite{CayleySubmanifoldsOfCalabiYauFourFolds}.

From the perspective of deformation theory, Ohst has examined boundary-value problems and generic transersality for compact Cayley submanifolds \cite{DeformationOfCompactCayleySubmanifoldsWithBoundary}, as well as the asymptotically cylindrical case \cite{DeformationsOfConicallySingularCayleySubmanifolds}. Further work has addressed deformations of compact complex surfaces regarded as Cayley submanifolds, and, by the same author, the conically singular setting \cite{CayleyDeformationsOfCompactComplexSurfaces}. More recently, Englebert has developed a broader programme studying the existence of calibrated fibrations in the $ \mathrm{Spin}(7) $ setting, together with the deformation theory of compact, conically singular, and asymptotically conical Cayley submanifolds, culminating in a generalised framework for the deformation theory of families of almost-Cayley submanifolds \cite{englebert2023conicallysingularcayleysubmanifoldsOne, englebert2023conicallysingularcayleysubmanifoldsTwo, englebert2024conicallysingularcayleysubmanifoldsThree}.

Cayley submanifolds also arise in theoretical physics. Becker, Becker, Morrison, Ooguri, Oz, and Yin show that Cayley cycles in $ \mathrm{Spin}(7) $-holonomy eight-manifolds, like associative and coassociative cycles in $ G_{2} $-holonomy seven-manifolds, preserve half of the spacetime supersymmetry \cite{becker1996supersymmetric}. By contrast, holomorphic and special Lagrangian cycles in Calabi-Yau $ 4 $-folds preserve half the spacetime supersymmetry as well, so what makes Cayley submanifolds distinctive is that they instead preserve only a quarter. The authors note that it would be interesting to find explicit examples of this phenomenon.

This paper therefore fits nicely into the existing literature, contributing a general method of construction, analysis, and understanding of singularities.

\subsection*{Characterisation of Cayley $ 4 $-Folds}

It useful to have a good understanding of a characterisation of Cayley $ 4 $-folds in $ \mathbb{C}^{4} $ for the work that follows.

\begin{lemma}[Harvey--Lawson, Lemma 4.4]
\label{LEMMA:4.4.Harvey-Lawson}

Given any $ 3 $-plane $ \eta $ in $ \mathbb{O} \cong \mathbb{R}^{8} $, there is a unique Cayley $ 4 $-plane $ \xi $ with $ \eta \subset \xi $.

\end{lemma}

\begin{theorem}

Suppose that $ N $ is a three-dimensional, real-analytic submanifold of $ \mathbb{O} $. Then there exists a unique real-analytic Cayley submanifold $ M $ of $ \mathbb{O} $ which contains $ N $ \cite[Theorem~4.5.]{harvey_calibrated_1982}

\end{theorem}

We now introduce a definition for the cross product of three Cayley numbers. 

\begin{definition}[The triple cross product of Cayley Numbers]\label{DEFINITION:TheTripleCrossProduct}

    Let $ (\mathbb{R}^{8}, g,\Phi) \cong (\mathbb{C}^{4}, g, \Phi) $ be equipped with the standard $ \mathrm{Spin}(7) $-structure $ \Phi=\frac{1}{2}\omega\wedge\omega+\operatorname{Re}\Omega $. The \textit{calibrated triple cross product} is the unique trilinear map $ \times:\Lambda^{3}\mathbb{R}^{8}\longrightarrow \mathbb{R}^{8} $ defined by the identity

    \begin{equation*}
        g(x\times y\times z,\xi)=\Phi(x,y,z,\xi), \quad \forall (x, y, z, \xi) \in \mathbb{R}^{8}.
    \end{equation*}

    Equivalently, if $ ({E_{\alpha}})_{\alpha=1}^{8} $ is any oriented orthonormal basis of $ \mathbb{R}^{8} $, then

    \begin{equation*}
        x \times y \times z = \sum_{\alpha = 1}^{8} \Phi(x, y, z, E_{\alpha})E_{\alpha}.
    \end{equation*}

    Under the identification $ \mathbb{R}^{8} \cong \mathbb{C}^{4} $, with real basis $ e_{1}, Je_{1}, \ldots, e_{4}, Je_{4} $, $ Je_{k} = ie_{k} $, the $ k $-th complex component is

    \begin{equation*}
        (x \times y \times z)_{k} = \Phi(x, y, z, e_{k}) + i \Phi(x, y, z, Je_{k}).
    \end{equation*}
    
\end{definition}

If one identifies $ \mathbb{R}^{8} $ with the octonions $ \mathbb{O} $, then this agrees with the octonionic formula $ x \times y \times z = \tfrac{1}{2} (x(\overline{y}z) - z(\overline{y}x)) $, as per \cite{harvey_calibrated_1982}.

This allows us to state and prove the following proposition.

\begin{proposition}\label{PROP:CayleyCrossProductProofMethod}

    Let $ W_{1}, W_{2}, W_{3} \in \mathbb{C}^{4} $ be linearly independent over $ \mathbb{R} $. Then $ W_{1}, W_{2}, W_{3}, W_{1} \times W_{2} \times W_{3} $ are linearly independent over $ \mathbb{R} $, and $ \langle W_{1}, W_{2}, W_{3}, W_{1} \times W_{2} \times W_{3} \rangle $ is the unique Cayley $ 4 $-plane containing $ \langle W_{1}, W_{2}, W_{3} \rangle $.
    
\end{proposition}

Under the identification $ \mathbb{R}^{8} \cong \mathbb{C}^{4} \cong \mathbb{O} $, our Cayley triple cross product agrees with that of Harvey-Lawson. The result therefore follows from  \cite[Chapter IV, Lemma~4.4.]{harvey_calibrated_1982}.

\section{Geometric Setting}\label{SECT:GeometricalSetting}

\subsection{The $ \mathrm{U}(1)^{2} $-Action}

We work with coordinates $ (z_{1}, z_{2}, z_{3}, z_{4}) \in \mathbb{C}^{4} \cong \mathbb{R}^{8} $, with $ \mathbb{C}^{4} $ equipped with its standard Calabi-Yau $ 4 $-fold structure $ (\omega, \Omega) $, where

\begin{equation*}
    \omega = \frac{i}{2} \sum_{j=1}^{4} dz_{j} \wedge d\overline{z}_{j}, \quad \Omega = dz_{1} \wedge dz_{2} \wedge dz_{3} \wedge dz_{4}.
\end{equation*}

The Cayley calibration is

\begin{equation}
    \Phi = \frac{1}{2} \omega \wedge \omega + \mathrm{Re}(\Omega),
\end{equation}

which is the standard $ \mathrm{Spin}(7) $-structure on $ \mathbb{C}^{4} $.

For $ (e^{i \theta_{1}}, e^{i \theta_{2}}) \in \mathrm{U(1)}^{2} $, we define the $ \mathrm{U(1)}^{2} $-action on $ \mathbb{C}^{4} $ by

\begin{equation}\label{EQ:MainAction}
    (e^{i\theta_{1}}, e^{i\theta_{2}}) \cdot (z_{1}, z_{2}, z_{3}, z_{4})
    = \bigl(e^{i\theta_{1}}z_{1},\; e^{i\theta_{2}}z_{2},\;
      e^{-i(\theta_{1}+\theta_{2})}z_{3},\; z_{4}\bigr).
\end{equation}

with moment map  

\begin{equation}\label{EQ:MomentMap}
    \mu(z) = \bigl(\lvert z_{1}\rvert^{2} - \lvert z_{3}\rvert^{2},\;
                   \lvert z_{2}\rvert^{2} - \lvert z_{3}\rvert^{2}\bigr).
\end{equation}

The $ \mathrm{U}(1)^{2} $-action \eqref{EQ:MainAction} preserves $ \omega $ as it acts unitarily, and the action preserves $ \Omega $ as the weights sum to zero. Thus, $ \Phi $ is invariant under this action.

\subsection{Invariant Coordinates and the Orbit Space}

We wish to describe $ \mathrm{U(1)}^{2} $-invariant submanifolds of $ \mathbb{C}^{4} $ via functions which descend to the quotient $ \mathbb{C}^{4} / \mathrm{U(1)}^{2} $. 

Observe that both $ z_{4} $ and $ iz_{1}z_{2}z_{3} $ are invariant under the action \eqref{EQ:MainAction}. We therefore choose the following complex representations of these invariants:

\begin{equation}\label{EQ:InvariantVariables}
    z_{4} = x + iu,\qquad iz_{1}z_{2}z_{3}=v+iy,
\end{equation}

where each of $ x, y, u, v $ are real-valued invariant functions. The choice of phase and the ordering of these variables are adapted to our subsequent description of invariant $ 4 $-folds as graphs over the $ (x, y) $-plane.

To describe the remaining invariant quantities, we write

\begin{equation*}
    z_{1} = \sqrt{w+a_{1}}e^{i \theta_{1}}, \quad z_{2} = \sqrt{w+a_{2}}e^{i \theta_{2}}, \quad z_{3} = \sqrt{w+a_{3}}e^{i \theta_{3}} \quad (\theta \in \mathbb{R}),
\end{equation*}

for which we can separate the invariants into their radial and moment-map components. Define $ a_{1}, a_{2}, a_{3}: \mathbb{C}^{4} \to \mathbb{R} $ by 

\begin{equation}
    \begin{aligned}
        a_{1} &= \frac{2}{3}|z_{1}|^{2} - \frac{1}{3}|z_{2}|^{2} - \frac{1}{3}|z_{3}|^{2}, \\
        a_{2} &= -\frac{1}{3}|z_{1}|^{2} + \frac{2}{3}|z_{2}|^{2} - \frac{1}{3}|z_{3}|^{2}, \\
        a_{3} &= -\frac{1}{3}|z_{1}|^{2} - \frac{1}{3}|z_{2}|^{2} + \frac{2}{3}|z_{3}|^{2},
    \end{aligned}
\end{equation}

such that the functions $ a_{1}, a_{2}, a_{3} $ satisfy the trace-free relation $ a_{1}+a_{2}+a_{3} =0 $, and we can recover the moment map by

\begin{equation}
    \begin{aligned}
        a_{1} - a_{3} &= |z_{1}|^{2} - |z_{3}|^{2}, \\
        a_{2} - a_{3} &= |z_{2}|^{2} - |z_{3}|^{2}.
    \end{aligned}
\end{equation}

The complementary radial variable is

\begin{equation}
    w = \frac{1}{3} \Bigl(|z_{1}|^{2} + |z_{2}|^{2} + |z_{3}|^{2} \Bigr),
\end{equation}

from which we recover

\begin{equation}\label{EQ:RaddiFromw}
    |z_{j}|^{2} = w + a_{j}, \qquad j = 1, 2, 3.
\end{equation}

Taking the real and imaginary parts of $ iz_{1}z_{2}z_{3} $ and $ z_{4} $ gives us the definitions for $ x, u, v, y $:

\begin{align*}
    x &= \mathrm{Re}(z_{4}), \quad u = \mathrm{Im}(z_{4}), \\
    v &= \mathrm{Re}(iz_{1}z_{2}z_{3}), \quad y = \mathrm{Im}(iz_{1}z_{2}z_{3}).
\end{align*}

Together with the two independent moment-map variables, we now have local coordinates $ (x, y, u, v, a_{1}, a_{2}) $ on the principal stratum of the orbit space $ \mathbb{C}^{4} / \mathrm{U(1)}^{2} $.

A $ \mathrm{U(1)}^{2} $-invariant $ 4 $-fold descends to a $ 2 $-dimensional surface in this orbit space. Define $ \Sigma \subset \mathbb{R}^{2} $, where $ (x, y) \in \Sigma $. We choose to restrict to $ N $ for which $ (x, y) $ are local coordinates on $ N/\mathrm{U(1)}^{2} $, writing that $ a_{1} = a_{1}(x, y) $, $ a_{2} = a_{2}(x, y) $, $ u = u(x, y) $, and $ v = v(x, y) $.

\subsection{Unique Solutions for $ w $ and Singularities}

The invariant variables introduced are not independent. Taking the squared modulus of $ v + iy = iz_{1}z_{2}z_{3} $, and using \eqref{EQ:RaddiFromw}, we obtain the invariant cubic equation

\begin{equation}\label{EQ:PDefinition}
    P(w) := (w+a_{1})(w+a_{2})(w+a_{3}) = v^{2} + y^{2}.
\end{equation}

We can thus realise $ w $ as a root of \eqref{EQ:PDefinition}. In particular, for fixed $ (a_{1}, a_{2}, a_{3}, v, y) $, a root $ w $ of \eqref{EQ:PDefinition} is called \textit{admissible} if 

\begin{equation}
    w \geq -\min_{j=1, 2, 3}(a_{j}),
\end{equation}

so that $ w + a_{j} \geq 0 $ for each $ j $.

On the admissible region, there exists a unique admissible solution $ w $ which depends smoothly on $ (v, y) $. Writing $ w_{0} = -\min_{j}(a_{j}) $, we have $ P(w_{0}) = -(v^{2}+y^{2}) \leq 0 $, while $ P(w) \to \infty $ as $ w \to \infty $. Moreover, 

\begin{equation}\label{EQ:PDerivative}
    P'(w)=
    \sum_{j=1}^3\prod_{k\neq j}(w+a_{k})>0
    \qquad\text{for }w>w_{0}.
\end{equation}

Thus, $ P $ is strictly increasing on the admissible half-line.

Whenever $ P'(w) > 0 $, the Implicit Function Theorem gives $ w $ as a smooth function of $ (a_{1}, a_{2}, v, y) $, with $ a_{3} = -a_{1}-a_{2} $. In particular, this holds whenever $ (v, y) \neq (0, 0) $,

It remains to analyse the locus $ (v, y) = (0, 0) $. At such points, the admissible root is $ w = -\min_{j}(a_{j}) $, and its regularity is governed by the multiplicity of the minimum. For a point $ p $, we therefore define

\begin{equation}\label{EQ:MultiplicityFormula}      m(p)=\#\bigl\{\, j\in\{1,2,3\} : a_j(p)=\min_k a_k(p) \,\bigr\},
\end{equation}

and distinguish the three cases.

\begin{itemize}

\item \textbf{$ m(p) = 1 $}: Without loss of generality, suppose that $ a_{1}(p) < a_{2}(p), a_{3}(p) $. The admissible root is then $ w = -a_{1}(p) $, at which $ w+a_{1} = 0 $, and $ P'(-a_{1}) = (a_{2}-a_{1})(a_{3} - a_{1}) > 0 $. Hence, the Implicit Function Theorem implies that $ w $ depends smoothly on the invariant variables near $ p $.

\item \textbf{$ m(p) = 2 $:} Suppose that $ a_{1}(p) = a_{2}(p) < a_{3}(p) $. The admissible root then satisfies $ w + a_{1} = w+a_{2} = 0 $, so $ P'(w) = 0 $ and the Implicit Function Theorem no longer applies. At these points, precisely two of $ z_{1}, z_{2}, z_{3} $ vanishes, which corresponds to the $ U(1)^{2} $-orbit degenerating to a $ U(1) $-orbit, and so the action has a one-dimensional stabiliser.

\item \textbf{$ m(p) = 3 $:} In the maximally degenerate case, $ a_{1}(p) = a_{2}(p) = a_{3}(p) = 0 $, and equation \eqref{EQ:PDefinition} reduces to $ w^{3} = v^{2} + y^{2} $, so $ w = (v^{2}+y^{2})^{1/3} $. Thus, the admissible-root function is not locally Lipschitz at $ p $. Geometrically, $ z_{1}=z_{2}=z_{3} = 0 $, and the entire $ U(1)^{2} $-action fixes a point.

\end{itemize}

Accordingly, we distinguish the regular and degenerate regimes as follows: the \textit{regular regime} consists of points with either $ P(w) > 0 $, or else $ v=y=0 $ and $ m(p) $ is at most one; the \textit{degenerate regime } occurs when $ v = y = 0 $ and $ m(p) $ is at least $ 2 $.

\subsection{The $ \mathrm{U(1)}^{2} $-Invariant Submanifold $ N$}

Let $ \Sigma \subset \mathbb{R}^{2} $ be an open domain, and let $ u, v, a_{1}, a_{2}: \Sigma \to \mathbb{R} $, with $ a_{3}=-a_{1}-a_{2} $. Using the invariant coordinates introduced above, a $ \mathrm{U}(1)^{2} $-invariant submanifold $ N \subset \mathbb{C}^{4} $ is realised as

\begin{equation}\label{EQ:NDef}
\begin{aligned}
N := \Bigl\{\, z \in \mathbb{C}^{4} \;\Big|\;
&|z_{1}|^{2} - |z_{3}|^{2}
= a_{1}(x,y) - a_{3}(x,y), \\
&|z_{2}|^{2} - |z_{3}|^{2}
= a_{2}(x,y) - a_{3}(x,y), \\
&iz_{1}z_{2}z_{3} = v(x, y) + iy, \\
&z_{4} = x + i\,u(x,y)
\Bigr\},
\end{aligned}
\end{equation}

where $ x = \mathrm{Re}(z_{4}) $ and $ y = \mathrm{Im}(iz_{1}z_{2}z_{3}) $. 

On the locus where the action is free, the invariant coordinates $ (x, y, u, v, a_{1}, a_{2}) $ describe the orbit space $ \mathbb{C}^{4} / U(1)^{2} $. Under these coordinates, $ N / U(1)^{2} $ is represented by the graph $ (x, y) \longmapsto \bigl(x, y, u(x, y), v(x, y), a_{1}(x, y), a_{2}(x, y) \bigr) $. Consequently, the free part of $ N $ is a smooth $ 4 $-dimensional submanifold of $ \mathbb{C}^{4} $, with $ (x, y) $ providing coordinates on its $ 2 $-dimensional quotient.

\section{A Construction of $ \mathrm{U(1)}^{2} $-Invariant Cayley $ 4 $-Folds in $ \mathbb{C}^{4} $}\label{SECT:ConstructionOfCayleyFourFolds}

In this section, we seek necessary and sufficient conditions on the functions $ u, v, a_{1}, a_{2} $ for the $ \mathrm{U}(1)^{2} $-invariant submanifold $ N $ to be calibrated by the Cayley form $ \Phi $.

\subsection{Statement of the Main Theorem}

We begin by stating the main result of this paper.

\begin{theorem}[A construction of  $ \mathrm{U(1)}^{2} $-invariant Cayley $ 4 $-folds]\label{THM:MainTheoremRegular}

Let $ \Sigma \subset \mathbb{R}^{2} $ be a smooth domain, and $ a_{1}, a_{2}, a_{3}, u, v: \Sigma \to \mathbb{R} $ be continuous functions which satisfy $ a_{1} + a_{2} + a_{3} = 0 $.

Define the degenerate set

\begin{equation*}
    S := \bigl\{ (x, y) \in \Sigma: y =v(x, y)=0, \quad m(p) \geq 2 \},
\end{equation*}

and $ \Sigma^{*} = \Sigma \backslash S $, so that on $ \Sigma^{*} $, we may assume that $ a_{1}, a_{2}, a_{3}, u, v $ are smooth. 

For each $ (x, y) \in \Sigma $, let $ w = w(x, y) $ denote the unique admissible root of $ P(w) = v^{2} + y^{2} = (w+a_{1})(w+a_{2})(w+a_{3}) $, where $ w \geq - \min_{j}(a_{j}(x, y)) $.
    
Define a smooth $ \mathrm{U(1)}^{2} $-invariant four-dimensional submanifold $ N \subset \mathbb{C}^{4} $ by

\begin{equation}\label{EQ:NSubmanifold}
\begin{aligned}
N := \Bigl\{\, \vec{z} \in \mathbb{C}^{4} \;\Big|\;
&|z_{1}|^{2} - |z_{3}|^{2}
= a_{1}(x,y) - a_{3}(x,y), \\
&|z_{2}|^{2} - |z_{3}|^{2}
= a_{2}(x,y) - a_{3}(x,y), \\
& iz_{1}z_{2}z_{3} = v(x, y) + iy, \\ 
& z_{4} = x + iu(x, y), \quad (x, y) \in \Sigma^{*}
\Bigr\}.
\end{aligned}
\end{equation}

The submanifold $ N $ is Cayley with respect to $ \Phi $ if and only if the PDE system

\begin{equation}\label{EQ:MainPDESystem}
   \begin{cases}
    u_{x} = v_{y}, \\
    v_{x} = -P'(w)u_{y}, \\
    a_{1,x} = -\bigl((a_{1}+a_{2}+w)a_{1,y}+(a_{1}+2w)a_{2,y}\bigr), \\
a_{2,x} = (a_{2}+2w)a_{1,y} + (a_{1}+a_{2}+w)a_{2,y},
\end{cases}
\end{equation}

holds on $ \Sigma^{*} $, where $ P'(w) $ denotes differentiation with respect to the polynomial variable $ w $. 

\end{theorem}

The proof of Theorem~\ref{THM:MainTheoremRegular} occupies the remainder of this section.

\subsection{The Gauge and the Tangent Frame}

Before we commence the proof, it is necessary to understand the gauge and the tangent frame with which we shall be working.

On each orbit we choose the representative for which $ z_{1}, z_{2}, z_{3} $ share a common phase, writing $ z_{j}=\sqrt{w+a_{j}}\,e^{i\theta} $ for $ j = 1, 2, 3 $. The gauge $ \theta $ is determined by

\begin{equation}\label{EQ:ThetaDef}
    e^{3 i \theta} = \frac{y-iv}{\sqrt{v^{2} + y^{2}}},
\end{equation}

where $ 3 \theta = \arg(y-iv) =  \arg(z_{1}z_{2}z_{3}) $, the position vector thus being

\begin{equation}\label{EQ:PositionVector}
    W(x,y) :=
    \Bigl(
      \sqrt{w+a_{1}}\,e^{i \theta},\;
      \sqrt{w+a_{2}}\,e^{i \theta},\;
      \sqrt{w+a_{3}}\,e^{i \theta},\;
      x + i\,u(x,y)
    \Bigr).
\end{equation}

The tangent space $ T_{p}N $ is spanned by four vectors: \textit{two orbital directions} which are generated by the infinitesimal $ \mathrm{U(1)}^{2} $-action (which we shall denote by $ W_{\phi_{1}} $ and $ W_{\phi_{2}} $), and \textit{two transverse directions} which are obtained by differentiating~\eqref{EQ:PositionVector} with respect to $ x $ and $ y $ (which we shall denote by $ W_{x} $ and $ W_{y} $).

\begin{lemma}[Tangent Vectors]\label{LEMMA:TangentVectors}

The tangent space $ T_{p}N $ at each $ p \in N $ is spanned by the following four tangent vectors:

\begin{align}
    W_{\phi_{1}} &= (iz_{1},\; 0,\; -iz_{3},\; 0),
    \qquad
    W_{\phi_{2}} = (0,\; iz_{2},\; -iz_{3},\; 0), \notag\\[6pt]
    W_{x} &= \Bigl(
        \Bigl(\frac{w_{x}+a_{1,x}}{2\sqrt{w+a_{1}}}
          + i\theta_{x}\sqrt{w+a_{1}}\Bigr)e^{i\theta},\;
         \Bigl(\frac{w_{x}+a_{2,x}}{2\sqrt{w+a_{2}}}
          + i\theta_{x}\sqrt{w+a_{2}}\Bigr)e^{i\theta}, \notag\\
        &\qquad
         \Bigl(\frac{w_{x}+a_{3,x}}{2\sqrt{w+a_{3}}}
          + i\theta_{x}\sqrt{w+a_{3}}\Bigr)e^{i\theta},\;
         1 + i\,u_{x}
    \Bigr), \notag\\[6pt]
    W_{y} &= \Bigl(
        \Bigl(\frac{w_{y}+a_{1,y}}{2\sqrt{w+a_{1}}}
          + i\theta_{y}\sqrt{w+a_{1}}\Bigr)e^{i\theta},\;
         \Bigl(\frac{w_{y}+a_{2,y}}{2\sqrt{w+a_{2}}}
          + i\theta_{y}\sqrt{w+a_{2}}\Bigr)e^{i\theta}, \notag\\
        &\qquad
         \Bigl(\frac{w_{y}+a_{3,y}}{2\sqrt{w+a_{3}}}
          + i\theta_{y}\sqrt{w+a_{3}}\Bigr)e^{i\theta},\;
         i\,u_{y}
    \Bigr)\label{eq:Wall}
\end{align}

where:

\begin{align}
    \theta_{x} &= -\frac{y\,v_{x}}{3(v^{2}+y^{2})},
    &\theta_{y} &= \frac{v-yv_{y}}{3(v^{2}+y^{2})},\notag\\[4pt]
    w_{x} &= \frac{2v\,v_{x}
              - \displaystyle\sum_{j=1}^{3} a_{j,x}\!\prod_{k \neq j}(w+a_{k})}
             {P'(w)},
    &w_{y} &= \frac{2v\,v_{y} + 2y
              - \displaystyle\sum_{j=1}^{3} a_{j,y}\!\prod_{k \neq j}(w+a_{k})}
             {P'(w)}.\label{eq:theta_w}
\end{align}

\end{lemma}

\begin{proof}

The one-parameter subgroups of the $ \mathrm{U(1)}^{2} $-action are 

\begin{equation*}
    \gamma_{1}(\phi) = (e^{i\phi}z_{1},\; z_{2},\; e^{-i\phi}z_{3},\; z_{4}),
    \qquad
    \gamma_{2}(\phi) = (z_{1},\; e^{i\phi}z_{2},\; e^{-i\phi}z_{3},\; z_{4}).
\end{equation*}

Differentiating at $ \phi = 0 $ gives the expression for $ W_{\phi_{1}} $ and $ W_{\phi_{2}} $. 

For the transverse vectors, we differentiate \eqref{EQ:PositionVector} with respect to $ x $ and $ y $, treating $ w $ and $ \theta$ as implicit functions. Differentiation of $ (w+a_{1})(w+a_{2})(w+a_{3}) = v^{2}+y^{2} $ with respect to $ x $ gives

\begin{equation*}
    P'(w)\,w_{x} = 2v\,v_{x} - \sum_{j=1}^3 a_{j,x}\prod_{k \neq j}(w+a_{k}),
\end{equation*}

with the formula for $ w_{x} $ following by division by $ P'(w) $. The formula for $ w_{y} $ follows similarly, with the additional term $ 2y $ arising from $ \partial_{y}(y^{2}) $. For the $ \theta $ derivatives, we use \eqref{EQ:ThetaDef} with $ 3 \theta = \arg(y-iv) $. Then,

\begin{align*}
    3\theta_{x} = \mathrm{Im}\Bigl( \frac{-iv_{x}}{y-iv} \Bigr) = - \frac{yv_{x}}{v^{2}+y^{2}}, \quad 3 \theta_{y} = \mathrm{Im}\Bigl( \frac{1-iv_{y}}{y-iv} \Bigr) = \frac{v-yv_{y}}{v^{2}+y^{2}.}
\end{align*}
    
\end{proof} 

\subsection{The Cayley Cross Product}

The proof of Theorem~\ref{THM:MainTheoremRegular} will be a more general version of the proof of Proposition 4.1. in \cite{joyce2005u1}, and the proof of Theorem 1 in \cite{MyPaper}. We therefore state the following proposition.

\begin{proposition}\label{PROPOSITION:TheCayleyTripleCrossProductWy}

    Let $ \vec{z} = (z_{1}, z_{2}, z_{3}, z_{4}) \in N $ be a non-singular point, and let $ W_{\phi_{1}}, W_{\phi_{2}}, W_{x}, W_{y} $ be linearly independent over $ \mathbb{R} $. Then $ T_{p}N = \langle W_{\phi_{1}}, W_{\phi_{2}}, W_{x}, W_{y} \rangle $, and $ N $ is Cayley at $ \vec{z} $ if and only if

    \begin{equation*}
    W_y \in \langle W_{\phi_{1}}, W_{\phi_{2}}, W_{x}, 
    W_{\phi_{1}} \times W_{\phi_{2}} \times W_{x} \rangle,
\end{equation*}

where $ \times: \mathbb{C}^{4} \times \mathbb{C}^{4} \times \mathbb{C}^{4} \to \mathbb{C}^{4} $ denotes the Cayley triple cross product as in Definition~\ref{DEFINITION:TheTripleCrossProduct}.
    
\end{proposition}

It is thus of immediate use to us to compute the cross product $ W_{\phi_{1}} \times W_{\phi_{2}} \times W_{x} $.

\begin{proposition}\label{PROPOSITION:ComputationoftheCayleyCrossProduct}

    Let $ W_{\phi_{1}}, W_{\phi_{2}}, W_{x} $ be defined as in Lemma~\ref{LEMMA:TangentVectors}, and let $ \times: \mathbb{C}^{4} \times \mathbb{C}^{4} \times \mathbb{C}^{4} \to \mathbb{C}^{4} $ be defined as in Definition~\ref{DEFINITION:TheTripleCrossProduct}. Then

\begin{align}
        W_{\phi_{1}} \times W_{\phi_{2}} \times W_{x} = \Bigl(e^{i \theta} \Bigl( \sqrt{w+a_{2}}\sqrt{w+a_{3}}(\frac{y+iv}{\sqrt{v^{2} + y^{2}}})(1-iu_{x}) + \frac{1}{2}\sqrt{w+a_{1}}(a_{2, x} - a_{3, x}) \Bigr) , \notag\\
        e^{i \theta} \Bigl( \sqrt{w+a_{1}}\sqrt{w+a_{3}}(\frac{y+iv}{\sqrt{v^{2} + y^{2}}})(1-iu_{x}) - \frac{1}{2}\sqrt{w+a_{2}}(a_{1, x} - a_{3, x}) \Bigr) , \notag\\
        e^{i \theta}\Bigl( \sqrt{w+a_{1}}\sqrt{w+a_{2}}(\frac{y+iv}{\sqrt{v^{2} + y^{2}}})(1-iu_{x}) + \frac{1}{2} \sqrt{w+a_{3}}(a_{1, x} - a_{2, x}), \Bigr) \notag\\
        -iv_{x}\Bigr)\label{EQ:CrossProductRestult}
 \end{align}
    
\end{proposition}

\begin{proof}

By Definition~\ref{DEFINITION:TheTripleCrossProduct}, under the identification $ \mathbb{R}^{8} \cong \mathbb{C}^{4} $, with real orthonormal basis $ e_{1}, Je_{1}, \dots, e_{4}, Je_{4} $ ($ Je_{j} = ie_{j} $), the $ j $-th complex component of the Cayley triple cross product is 

\begin{equation}\label{EQ}
\begin{aligned}
\bigl(W_{\phi_{1}}\times W_{\phi_{2}}\times W_{x}\bigr)_{j}
={}&
\Phi(W{\phi_{1}},W_{\phi_{2}},W_{x},e_{j})
+i\Phi(W_{\phi_{1}},W_{\phi_{2}},W_{x},ie_{j}).
\end{aligned}
\end{equation}

Since $ \Phi=\frac12\omega\wedge\omega+\operatorname{Re}\Omega $, we compute separately the contributions from $ \mathrm{Re}(\Omega) $ and $ \tfrac{1}{2} \omega \wedge \omega $.

\medskip
\noindent\textbf{Calculation of $\mathrm{Re}(\Omega)$.}

For $ k = 1, 2, 3, 4 $, let $ M_{j} $ denote the $ 4 \times 4 $ complex matrices with columns $ W_{\phi_{1}}, W_{\phi_{2}}, W_{x}, e_{j} $, and define $ D_{j} := \det(M_{j}) = \Omega(W_{\phi_{1}}, W_{\phi_{2}}, W_{x}, e_{j}) $. Since $ \Omega $ is complex multilinear, $ \Omega(W_{\phi_1},W_{\phi_2},W_x,ie_j)=iD_j $, and hence 

\begin{equation*}
    \operatorname{Re}(D_j)
+i\,\operatorname{Re}(iD_j)
=
\overline{D_j}.
\end{equation*}

The complex vector representing the $ \mathrm{Re}(\Omega) $-contribution to the triple product is thus

\begin{equation*}
    \bigl(\overline{\det(M_1)}, \overline{\det(M_2)}, \overline{\det(M_3)}, \overline{\det(M_4)} \bigr).
\end{equation*}

Using cofactor expansion, the determinants are:

\begin{equation*}
\begin{aligned}
    \det(M_{1}) &= \sqrt{w+a_{2}}\sqrt{w+a_{3}}e^{2 i \theta}(1+iu_{x}), \\
    \det(M_{2}) &= \sqrt{w+a_{1}}\sqrt{w+a_{3}}e^{2i\theta}(1+iu_{x}), \\
    \det(M_{3}) &= \sqrt{w+a_{1}}\sqrt{w+a_{2}}e^{2 i \theta}(1+iu_{x}), \\
    \det(M_{4}) &= e^{3i\theta} \Bigl(- \tfrac{1}{2\sqrt{w+a_{3}}} \sqrt{w+a_{1}}\sqrt{w+a_{2}}(a_{3, x} + w_{x}) - \tfrac{1}{2\sqrt{w+a_{2}}} \sqrt{w+a_{1}}\sqrt{w+a_{3}}(a_{2, x}+w_{x}) \\
    &- \tfrac{1}{2\sqrt{w+a_{1}}} \sqrt{w+a_{2}}\sqrt{w+a_{3}}(a_{1, x}+w_{x})  -3i\theta_{x}\sqrt{w+a_{1}}\sqrt{w+a_{2}}\sqrt{w+a_{3}} \Bigr).
\end{aligned}
\end{equation*}

We are able to simplify the expression for $ \det(M_{4}) $ by noting that $ \det(M_{4}) $ is $ -(W_{x})_{1}z_{2}z_{3} - z_{1}(W_{x})_{2}z_{3} - z_{1}z_{2}(W_{x})_{3} $, and since $ (W_{x})_{j} $ is precisely the $ x $-derivative of $ z_{j} $ along the chosen parameterisation, we have that:

\begin{equation*}
    \det(M_{4}) = - \frac{\partial}{\partial x}(z_{1}z_{2}z_{3}).
\end{equation*}

Consequently, $ \det(M_{4}) = iv_{x} $.

Taking the complex conjugate of each of $ \det(M_{1}) $, $ \det(M_{2}) $, $ \det(M_{3}) $, and $ \det(M_{4}) $, one reaches the claimed expression for $ \mathrm{Re}(\Omega) $.

\medskip
\noindent\textbf{Calculation of $ \tfrac{1}{2}\omega\wedge\omega $.}

The term $ \tfrac{1}{2} \omega \wedge \omega $ can be expanded:

\begin{equation*}
    \frac{1}{2} \omega \wedge \omega(W_{\phi_{1}}, W_{\phi_{2}}, W_{x}, e_{i}) = \omega(W_{\phi_{1}}, W_{\phi_{2}})\omega(W_{x}, e_{i}) - \omega(W_{\phi_{1}}, W_{x})\omega(W_{\phi_{2}}, e_{i}) + \omega(W_{\phi_{2}}, W_{x})\omega(W_{\phi_{1}}, e_{i}).
\end{equation*}

Therefore, the computation for $ \tfrac{1}{2} \omega \wedge \omega $ reduces to computing each of the above terms.

We begin by noting that $ \omega(W_{\phi_{1}}, W_{\phi_{2}}) = 0 $, and so the first term vanishes. As for the other terms, direct computation yields

\begin{equation*}
    \omega(W_{\phi_{1}}, W_{x}) =-  \frac{1}{2}(a_{1, x} - a_{3, x}), \qquad \omega(W_{\phi_{2}}, W_{x}) = - \frac{1}{2}(a_{2, x} - a_{3, x}).
\end{equation*}

Therefore,

\begin{equation*}
    \tfrac{1}{2} \omega \wedge \omega(W_{\phi_{1}}, W_{\phi_{2}}, W_{x}, e_{i}) = \tfrac{1}{2}(a_{1, x} - a_{3, x})\omega(W_{\phi_{2}}, e_{i}) - \tfrac{1}{2}(a_{2, x} - a_{3, x})\omega(W_{\phi_{1}}, e_{i}).
\end{equation*}

Computing each of $ \omega(W_{\phi_{i}}, e_{i}) $ gives: 

\begin{align*}
  \omega(W_{\phi_{1}}, e_{1}) + i \omega(W_{\phi_{1}}, ie_{1}) = -z_{1}, \quad \omega(W_{\phi_{1}}, e_{2}) + i\omega(W_{\phi_{1}}, ie_{2}) = 0, \\\omega(W_{\phi_{1}}, e_{3}) + i\omega(W_{\phi_{1}}, ie_{3})= z_{3}, \quad \omega(W_{\phi_{1}}, e_{4}) + i\omega(W_{\phi_{1}}, ie_{4}) = 0,
\end{align*}

\begin{align*}
  \omega(W_{\phi_{2}}, e_{1}) + i \omega(W_{\phi_{1}}, ie_{1}) &= 0, \quad \omega(W_{\phi_{2}}, e_{2}) + i \omega(W_{\phi_{2}}, ie_{2}) = -z_{2}, \\ \omega(W_{\phi_{2}}, e_{3}) + i \omega(W_{\phi_{2}}, ie_{3}) &= z_{3}, \quad \omega(W_{\phi_{2}}, e_{4}) + i \omega(W_{\phi_{2}}, ie_{4}) = 0.
\end{align*}

Therefore, the four components of $ \tfrac{1}{2} \omega \wedge \omega $ are:

\begin{equation*}
    \begin{aligned}
        \frac{1}{2} \omega \wedge \omega(W_{\phi_{1}}, W_{\phi_{2}}, W_{x}, e_{1}) + i \frac{1}{2} \omega \wedge \omega (W_{\phi_{1}}, W_{\phi_{2}}, W_{x}, e_{1})  &=  \frac{1}{2}z_{1}(a_{2, x} - a_{3, x}), \\
        \frac{1}{2} \omega \wedge \omega(W_{\phi_{1}}, W_{\phi_{2}}, W_{x}, e_{2}) + i \frac{1}{2} \omega \wedge \omega (W_{\phi_{1}}, W_{\phi_{2}}, W_{x}, e_{2})   &= - \frac{1}{2}z_{2}(a_{1, x} - a_{3, x}), \\
        \frac{1}{2} \omega \wedge \omega(W_{\phi_{1}}, W_{\phi_{2}}, W_{x}, e_{3}) + i \frac{1}{2}  \omega \wedge \omega (W_{\phi_{1}}, W_{\phi_{2}}, W_{x}, e_{3})  &= \frac{1}{2}z_{3}(a_{1, x} - a_{2, x}), \\
        \frac{1}{2} \omega \wedge \omega(W_{\phi_{1}}, W_{\phi_{2}}, W_{x}, e_{4}) + i \frac{1}{2} \omega \wedge \omega(W_{\phi_{1}}, W_{\phi_{2}}, W_{x}, e_{4}) &= 0,
    \end{aligned}
\end{equation*}

with resulting complex vector

\begin{equation*}
    \Bigl( \frac{1}{2}z_{1}(a_{2, x} - a_{3, x}), -\frac{1}{2}z_{2}(a_{1, x} - a_{3, x}), \frac{1}{2}z_{3}(a_{1, x} - a_{2, x}), 0 \Bigr).
\end{equation*}

This gives the claimed form.

\end{proof}

\subsection{Proof of Theorem~\ref{THM:MainTheoremRegular}}

We are now in a position to prove Theorem~\ref{THM:MainTheoremRegular}. 

\begin{proof}

\textbf{Proof that the PDEs are necessary.}

Suppose that $ N $ is Cayley. By Proposition~\ref{PROP:CayleyCrossProductProofMethod}, there thus exists $ \alpha, \beta, \gamma, \delta \in \mathbb{R} $ such that

\begin{equation}\label{EQ:CayleyProofEQ}
    W_{y} = \alpha W_{\phi_{1}} + \beta W_{\phi_{2}} + \gamma W_{x} + \delta\, W_{\phi_{1}} \times W_{\phi_{2}} \times W_{x}.
\end{equation}

The task is therefore to compute each of $ \alpha, \beta, \gamma, \delta $, before substituting each back into \eqref{EQ:CayleyProofEQ} to derive the four PDEs.

\textit{Computation of $ \gamma $.} 

Equating the fourth component of each tangent vector in \eqref{EQ:CayleyProofEQ} gives  $ iu_{y} = \gamma(1+iu_{x}) -i \delta v_{x} $, for which taking the real parts determines that $ \gamma = 0 $, and taking the imaginary part determines that $ u_{y} = -\delta v_{x} $.

\textit{Computation of $ \delta $.} 

Equating the first three components of each tangent vector in \eqref{EQ:CayleyProofEQ}, we arrive at the equations:

\begin{subequations}\label{EQ:CayleyComponents}
\begin{align}
\begin{split}
    \frac{(w_{y} + a_{1, y})e^{i \theta}}{2 \sqrt{w + a_{1}}} + i \theta_{y} \sqrt{w + a_{1}} e^{i \theta}
        &= i \alpha \sqrt{w+a_{1}}e^{i \theta} \\
        &\quad + \delta\Bigl(\sqrt{w + a_{2}} \sqrt{w + a_{3}}\, e^{-2i \theta} (1-iu_{x}) + \tfrac{1}{2} \sqrt{w + a_{1}}\,e^{i \theta}(a_{2, x} - a_{3, x})\Bigr)
\end{split}
\label{EQ:CayleyComponentOne}\\[4pt]
\begin{split}
    \frac{(w_{y} + a_{2, y})e^{i \theta}}{2 \sqrt{w + a_{2}}} + i \theta_{y} \sqrt{w + a_{2}} e^{i \theta}
        &= i \beta \sqrt{w+a_{2}}e^{i \theta} \\
        &\quad + \delta\Bigl(\sqrt{w + a_{1}} \sqrt{w + a_{3}}\, e^{-2i \theta} (1-iu_{x}) - \tfrac{1}{2} \sqrt{w + a_{2}}\,e^{i \theta}(a_{1, x} - a_{3, x})\Bigr)
\end{split}
\label{EQ:CayleyComponentTwo}\\[4pt]
\begin{split}
    \frac{(w_{y} + a_{3, y})e^{i \theta}}{2 \sqrt{w + a_{3}}} + i \theta_{y} \sqrt{w + a_{3}} e^{i \theta}
        &= -i \sqrt{w+a_{3}}e^{i \theta}(\alpha + \beta) \\
        &\quad + \delta\Bigl(\sqrt{w + a_{1}} \sqrt{w + a_{2}}\, e^{-2i \theta} (1-iu_{x}) + \tfrac{1}{2} \sqrt{w + a_{3}}\,e^{i \theta}(a_{1, x} - a_{2, x})\Bigr)
\end{split}
\label{EQ:CayleyComponentThree}
\end{align}
\end{subequations}

This admits a normalisation: divide \eqref{EQ:CayleyComponentOne} by $ z_{1} $, divide \eqref{EQ:CayleyComponentTwo} by $ z_{2} $, and divide \eqref{EQ:CayleyComponentThree} by $ z_{3} $, before applying the formula for $ e^{-3i\theta} $ as per \eqref{EQ:ThetaDef}. The three normalised equations are then:

\begin{subequations}\label{EQ:NormalisedCayleyComponents}
\begin{align}
\frac{w_y+a_{1,y}}{2(w+a_1)}+i\theta_y
&=
i\alpha
+\delta\left(
\frac{(y+iv)(1-iu_x)}{w+a_1}
+\frac{a_{2,x}-a_{3,x}}{2}
\right),
\label{EQ:NormalisedCayleyComponentOne}\\[4pt]
\frac{w_y+a_{2,y}}{2(w+a_2)}+i\theta_y
&=
i\beta
+\delta\left(
\frac{(y+iv)(1-iu_x)}{w+a_2}
-\frac{a_{1,x}-a_{3,x}}{2}
\right),
\label{EQ:NormalisedCayleyComponentTwo}\\[4pt]
\frac{w_y+a_{3,y}}{2(w+a_3)}+i\theta_y
&=
-i(\alpha+\beta)
+\delta\left(
\frac{(y+iv)(1-iu_x)}{w+a_3}
+\frac{a_{1,x}-a_{2,x}}{2}
\right).
\label{EQ:NormalisedCayleyComponentThree}
\end{align}
\end{subequations}

By summing the equations~\eqref{EQ:NormalisedCayleyComponents},
the terms involving $ \alpha $ and $ \beta $, as well as the differences $ a_{j,x}-a_{k,x} $, cancel. This leaves:

\begin{equation}\label{EQ:ProofIntermediateStepOne}
\begin{aligned}
    \frac{w_{y} + a_{1, y}}{2(w+a_{1})} + \frac{w_{y} + a_{2, y}}{2(w+a_{2})} + \frac{w_{y} + a_{3, y}}{2(w+a_{3})} + 3i \theta_{y} = \\\delta \Bigl( \frac{(y+iv)(1-iu_{x})}{w+a_{1}} + \frac{(y+iv)(1-iu_{x})}{w+a_{2}} + \frac{(y+iv)(1-iu_{x})}{w+a_{3}} \Bigr).
\end{aligned}
\end{equation}

The left-hand-side of \eqref{EQ:ProofIntermediateStepOne} is equivalent to $ (z_{1}z_{2}z_{3})_{y}/z_{1}z_{2}z_{3} = (1-iv_{y})/(y-iv) $. The right-hand-side of \eqref{EQ:ProofIntermediateStepOne} is equivalent to $ \delta P'(w) (y+iv)(1-iu_{x})/P(w) $. Therefore, we are left with $ 1-iv_{y} = \delta P'(w)(1-iu_{x}) $, for which taking the real parts gives us that $ \delta = 1/P'(w) $.

\textit{Computation of $ \alpha $ and $ \beta $.} 

Consider the first normalised equation \eqref{EQ:NormalisedCayleyComponentOne}. Substituting in $ \delta = 1/P'(w) $ gives:  

\begin{equation}\label{EQ:ProofIntermediateStepThree}
    \frac{(w_{y} + a_{1, y})}{2(w+a_{1})} + i\theta_{y} = i \alpha + \frac{1}{P'(w)} \Bigl( \frac{(y+iv)(1-iu_{x})}{(w+a_{1})} + \frac{a_{2, x} - a_{3, x}}{2} \Bigr).
\end{equation}

Taking the imaginary parts of \eqref{EQ:ProofIntermediateStepThree}, and then taking the imaginary parts of \eqref{EQ:NormalisedCayleyComponentTwo}, determines $ \alpha $ and $ \beta $:

\begin{equation*}
    \alpha = \theta_{y} - \frac{v-yu_{x}}{P'(w)(w+a_{1})}, \quad \beta = \theta_{y} - \frac{v-yu_{x}}{P'(w)(w+a_{2})}.
\end{equation*}

This completes the calculation of $ \alpha, \beta, \gamma, \delta $.

\textit{Derivation of the reduced system.} 

Now substitute each of $ \alpha, \beta, \gamma, \delta $ into \eqref{EQ:NormalisedCayleyComponents}, leaving

\begin{subequations}\label{EQ:SubstitutedCayleyComponents}
\begin{align}
\frac{w_y+a_{1,y}}{2(w+a_1)}+i\theta_y
&=
i\left(
\theta_y-\frac{v-yu_x}{P'(w)(w+a_1)}
\right) \notag\\
&\quad
+\frac{1}{P'(w)}
\left(
\frac{(y+iv)(1-iu_x)}{w+a_1}
+\frac{a_{2,x}-a_{3,x}}{2}
\right),
\label{EQ:SubstitutedCayleyComponentOne}\\[6pt]
\frac{w_y+a_{2,y}}{2(w+a_2)}+i\theta_y
&=
i\left(
\theta_y-\frac{v-yu_x}{P'(w)(w+a_2)}
\right) \notag\\
&\quad
+\frac{1}{P'(w)}
\left(
\frac{(y+iv)(1-iu_x)}{w+a_2}
-\frac{a_{1,x}-a_{3,x}}{2}
\right),
\label{EQ:SubstitutedCayleyComponentTwo}\\[6pt]
\frac{w_y+a_{3,y}}{2(w+a_3)}+i\theta_y
&=
-i\left(
2\theta_y
-\frac{v-yu_x}{P'(w)(w+a_1)}
-\frac{v-yu_x}{P'(w)(w+a_2)}
\right) \notag\\
&\quad
+\frac{1}{P'(w)}
\left(
\frac{(y+iv)(1-iu_x)}{w+a_3}
+\frac{a_{1,x}-a_{2,x}}{2}
\right), \label{EQ:SubstitutedCayleyComponentThree}\\
iu_y=-\frac{i v_x}{P'(w)}\label{EQ:SubstitutedCayleyComponentFour}
\end{align}
\end{subequations}

Separating the real and imaginary parts of \eqref{EQ:SubstitutedCayleyComponents} gives us a system of eight equations, of which, three are tautologies:

\begin{subequations}\label{EQ:ReducedComponentEquations}
\begin{align}
\frac{w_y+a_{1,y}}{2(w+a_1)}
&=
\frac{1}{P'(w)}
\left(
\frac{y+vu_x}{w+a_1}
+\frac{a_{2,x}-a_{3,x}}{2}
\right),
\label{EQ:ReducedRealComponentOne}
\\[4pt]
\theta_y
&=
\theta_y
-\frac{v-yu_x}{P'(w)(w+a_1)}
+\frac{v-yu_x}{P'(w)(w+a_1)} \quad \implies 0 = 0,
\label{EQ:ReducedImaginaryComponentOne}
\\[4pt]
\frac{w_y+a_{2,y}}{2(w+a_2)}
&=
\frac{1}{P'(w)}
\left(
\frac{y+vu_x}{w+a_2}
-\frac{a_{1,x}-a_{3,x}}{2}
\right),
\label{EQ:ReducedRealComponentTwo}
\\[4pt]
\theta_y
&=
\theta_y
-\frac{v-yu_x}{P'(w)(w+a_2)}
+\frac{v-yu_x}{P'(w)(w+a_2)} \quad \implies 0 = 0,
\label{EQ:ReducedImaginaryComponentTwo}
\\[4pt]
\frac{w_y+a_{3,y}}{2(w+a_3)}
&=
\frac{1}{P'(w)}
\left(
\frac{y+vu_x}{w+a_3}
+\frac{a_{1,x}-a_{2,x}}{2}
\right),
\label{EQ:ReducedRealComponentThree}
\\[4pt]
\theta_y
&=
-2\theta_y
+\frac{v-yu_x}{P'(w)(w+a_1)}
+\frac{v-yu_x}{P'(w)(w+a_2)}
+\frac{v-yu_x}{P'(w)(w+a_3)},
\label{EQ:ReducedImaginaryComponentThree}
\\[4pt]
u_y
&=
-\frac{1}{P'(w)}v_x.
\label{EQ:ReducedFourthComponent}\\
0&=0
\end{align}
\end{subequations}

The PDE $ u_{y} = - \tfrac{1}{P'(w)}v_{x} $ is immediate from \eqref{EQ:ReducedFourthComponent}. 

Considering \eqref{EQ:ReducedImaginaryComponentThree}, we combine the $ \theta_{y} $ terms on the left-hand-side to give

\begin{equation*}
    3\theta_{y} = \frac{v-yu_{x}}{P'(w)(w+a_{1})} + \frac{v-yu_{x}}{P'(w)(w+a_{2})} + \frac{v-yu_{x}}{P'(w)(w+a_{3})} = \frac{P'(w)(v-yu_{x})}{P'(w)(w+a_{1})(w+a_{2})(w+a_{3})}.
\end{equation*}

Therefore, 

\begin{equation*}
    3\theta_{y} = \frac{v-yv_{y}}{P(w)}.
\end{equation*}

Substituting in \eqref{EQ:ThetaDef}, we determine that $ y(u_{x} - v_{y}) = 0 $, and hence $ u_{x} = v_{y} $ on $ y \neq 0 $, which extends across $ y = 0 $ by continuity.

Since $ u_{x} = v_{y} $, equation \eqref{EQ:ReducedRealComponentOne} can be compressed using that $ 2(y+vu_{x}) = 2(y+vv_{y}) $. By Lemma~\ref{LEMMA:TangentVectors}, $ 2(y+vu_{x}) = P'(w)w_{y} + (w+a_{2})(w+a_{3})a_{1, y} + (w+a_{1})(w+a_{3})a_{2, y} + (w+a_{1})(w+a_{2})a_{3, y} $. Thus, substituting this expression for $ 2(y+vu_{x}) $ into \eqref{EQ:ReducedRealComponentOne}--\eqref{EQ:ReducedRealComponentThree} cancels the $ w_{y} $ terms, leaving a linear system for the $ a_{i, x} $,

\begin{align*}
    P'(w)a_{1, y} = (w+a_{2})(w+a_{3})a_{1, y}+(w+a_{1})(w+a_{3})a_{2, y} + (w+a_{1})(w+a_{2})a_{3, y} + (w+a_{1})(a_{2, x} - a_{3, x}), \\
    P'(w)a_{2, y} = (w+a_{2})(w+a_{3})a_{1, y} + (w+a_{1})(w+a_{3})a_{2, y} + (w+a_{1})(w+a_{2})a_{3, y} - (w+a_{2})(a_{1, x} - a_{3, x}), \\
    P'(w)a_{3, y} = (w+a_{2})(w+a_{3})a_{1, y} + (w+a_{1})(w+a_{3})a_{2, y} + (w+a_{1})(w+a_{2})a_{3, y} + (w+a_{3})(a_{1, x} - a_{2, x}).
\end{align*}

With the trace-free relation $ a_{1}+a_{2}+a_{3} = 0 $, solving the resultant linear system determines that

\begin{align*}
    a_{1, x} &= -\Bigl( (a_{1}+a_{2}+w)a_{1, y} + (a_{1}+2w)a_{2, y} \Bigr), \\
    a_{2, x} &= (a_{2}+2w)a_{1, y} + (a_{1}+a_{2}+w)a_{2, y}.
\end{align*}

This gives us the PDE system claimed in \eqref{EQ:MainPDESystem}, and we are able to conclude that if $ N $ is Cayley, the PDE system \eqref{EQ:MainPDESystem} must hold.

\noindent\textbf{Proof that the PDEs are sufficient.}

Assume that the PDE system \eqref{EQ:MainPDESystem} holds, that is, that the functions $ u, v, a_{1}, a_{2} $ are solutions of \eqref{EQ:MainPDESystem}. Using these four equations, we wish to recover the reduced component equations obtained from the Cayley cross-product identity. 

Define

\begin{equation*}
    \alpha = \theta_{y} - \frac{v-yu_{x}}{P'(w)(w+a_{1})}, \quad \beta = \theta_{y} - \frac{v-yu_{x}}{P'(w)(w+a_{2})}, \quad
    \gamma = 0, \quad \delta = \frac{1}{P'(w)}.
\end{equation*}

Using the definition of $ \theta_{y} $, these equations become \eqref{EQ:ReducedImaginaryComponentOne} and \eqref{EQ:ReducedImaginaryComponentTwo}, respectively.

Since $ u_{x} = v_{y} $, the formula for $ \theta_{y} $  becomes $ \theta_{y} = (v-yv_{y})/3(v^{2}+y^{2}) $, which is equivalent to $ \theta_{y} = (v-yu_{x})/3(v^{2}+y^{2}) $. The reciprocal identity

\begin{equation*}
    \frac{1}{w+a_{1}} + \frac{1}{w+a_{2}} + \frac{1}{w+a_{3}} = \frac{P'(w)}{(w+a_{1})(w+a_{2})(w+a_{3})}
\end{equation*} 

means that

\begin{equation*}
    3\theta_y
=
\frac{v-yu_x}{P'(w)}
\left(
\frac{1}{w+a_1}
+\frac{1}{w+a_2}
+\frac{1}{w+a_3}
\right),
\end{equation*}

or equivalently,

\begin{equation*}
    \theta_y
=
-2\theta_y
+\frac{v-yu_x}{P'(w)(w+a_1)}
+\frac{v-yu_x}{P'(w)(w+a_2)}
+\frac{v-yu_x}{P'(w)(w+a_3)}.
\end{equation*}

This determines \eqref{EQ:ReducedImaginaryComponentThree}.

It remains to verify that the three reduced real component equations imply

\begin{align*}
    P'(w)a_{1, y} = (w+a_{2})(w+a_{3})a_{1, y}+(w+a_{1})(w+a_{3})a_{2, y} + (w+a_{1})(w+a_{2})a_{3, y} + (w+a_{1})(a_{2, x} - a_{3, x}), \\
    P'(w)a_{2, y} = (w+a_{2})(w+a_{3})a_{1, y} + (w+a_{1})(w+a_{3})a_{2, y} + (w+a_{1})(w+a_{2})a_{3, y} - (w+a_{2})(a_{1, x} - a_{3, x}), \\
    P'(w)a_{3, y} = (w+a_{2})(w+a_{3})a_{1, y} + (w+a_{1})(w+a_{3})a_{2, y} + (w+a_{1})(w+a_{2})a_{3, y} + (w+a_{3})(a_{1, x} - a_{2, x}).
\end{align*}

By the PDEs governing $ a_{1} $ and $ a_{2} $ as per \eqref{EQ:MainPDESystem}, we may use that $ a_{3, x} = -a_{1, x} - a_{2, x} $, to determine that

\begin{align*}
   a_{2, x} - a_{3, x} &= (-a_{1}+a_{2}+3w)a_{1, y} + (a_{1}+2a_{2})a_{2, y}, \\
   a_{1, x} - a_{3, x} &= -(2a_{1}+a_{2})a_{1, y} + (-a_{1}+a_{2}-3w)a_{2, y}, \\
   a_{1, x}-a_{2, x} &= -(a_{1}+2a_{2}+3w)a_{1, y} - (2a_{1}+a_{2}+3w)a_{2, y}
\end{align*}

Furthermore,

\begin{align*}
    (w+a_{2})(w+a_{3})a_{1, y} + (w+a_{1})(w+a_{3})a_{2, y} + (w+a_{1})(w+a_{2})a_{3, y} \\
    = -(w+a_{2})(2a_{1}+a_{2})a_{1, y} - (w+a_{1})(a_{1}+2a_{2})a_{2, y},
\end{align*}

and 

\begin{align*}
    P'(w) = 3w^{2}-a_{1}^{2}- a_{1}a_{2}-a_{2}^{2}.
\end{align*}

Substitution gives

\begin{align}
\label{EQ:IntermediateProofStepFour}
    &(w+a_{2})(w+a_{3})a_{1,y}
    +(w+a_{1})(w+a_{3})a_{2,y} \notag\\
    &\quad +(w+a_{1})(w+a_{2})a_{3,y}
    +(w+a_{1})(a_{2,x}-a_{3,x})
    =P'(w)a_{1,y},\\[6pt]
    &(w+a_{2})(w+a_{3})a_{1,y}
    +(w+a_{1})(w+a_{3})a_{2,y} \notag\\
    &\quad +(w+a_{1})(w+a_{2})a_{3,y}
    +(w+a_{3})(a_{1,x}-a_{2,x})
    =P'(w)a_{3,y}.
\end{align}

Since $ u_{x} = v_{y} $, we can rewrite the formula for $ w_{y} $ from Lemma~\ref{LEMMA:TangentVectors} as

\begin{align*}
    w_{y}P'(w) = 2vu_{x} + 2y - a_{1, y}(w+a_{2})(w+a_{3}) - a_{2, y}(w+a_{1})(w+a_{3}) - a_{3, y}(w+a_{1})(w+a_{2}).
\end{align*}

Combining this with \eqref{EQ:IntermediateProofStepFour},

\begin{align*}
    P'(w)a_{1, y} = 2vu_{x}+2y - w_{y}P'(w) + (w+a_{1})(a_{2, x}-a_{3, x}) \\
    \implies w_{y}P'(w) + P'(w)a_{1, y} = 2vu_{x}+2y + (w+a_{1})(a_{2, x}-a_{3, x}) \\
    \frac{w_{y}+a_{1, y}}{2(w+a_{1})} = \frac{1}{P'(w)} \Bigl( \frac{y+vu_{x}}{w+a_{1}} + \frac{a_{2, x} - a_{3, x}}{2} \Bigr).
\end{align*}

Similarly,

\begin{align*}
    \frac{w_y+a_{2,y}}{2(w+a_2)}
&=
\frac{1}{P'(w)}
\left(
\frac{y+vu_x}{w+a_2}
-\frac{a_{1,x}-a_{3,x}}{2}
\right), \\
\frac{w_y+a_{3,y}}{2(w+a_3)}
&=
\frac{1}{P'(w)}
\left(
\frac{y+vu_x}{w+a_3}
+\frac{a_{1,x}-a_{2,x}}{2}
\right).
\end{align*}

These are the three reduced real component equations \eqref{EQ:ReducedRealComponentOne}--\eqref{EQ:ReducedRealComponentThree}. 

We therefore have verified all seven equations in \eqref{EQ:ReducedComponentEquations}. By the definitions of $ \alpha, \beta, \gamma, \delta $, these equations are precisely the real and imaginary components of $ W_y
=
\alpha W_{\phi_1}
+\beta W_{\phi_2}
+\gamma W_x
+\delta\,
W_{\phi_1}\times W_{\phi_2}\times W_x $. Moreover, $ \delta = \tfrac{1}{P'(w)} > 0 $. It follows from Proposition~\ref{PROP:CayleyCrossProductProofMethod} that $ N $ is Cayley on the locus $ v^{2}+y^{2} > 0 $, with the prescribed orientation. The Cayley identity extends to all of $ N $ by continuity. Hence, the PDE system \eqref{EQ:MainPDESystem} implies that $ N $ is Cayley.
\end{proof}

\section{Analysis of $ \mathrm{U(1)}^{2} $-Invariant Cayley $ 4 $-Folds in $ \mathbb{R}^{8} $}\label{SECT:AnalysisOfCayley4FoldConstruction}

\subsection{The Nature of the Cayley PDE System}

The PDE system \eqref{EQ:MainPDESystem} is a non-linear, first-order Cauchy-Riemann system. We are able to prove the ellipticity of the system.

\begin{proposition}[Ellipticity of the PDE system \eqref{EQ:MainPDESystem}]

    The PDE system \eqref{EQ:MainPDESystem} is elliptic precisely where $ P'(w) > 0 $, with the ellipticity degenerating as $ P'(w) \to 0 $.
    
\end{proposition}

\begin{proof}

    The $ (u, v) $-subsystem of \eqref{EQ:MainPDESystem} has the coefficient matrix

    \begin{equation*}
        A = \begin{pmatrix}
            0 & 1 \\
            -P'(w) & 0
        \end{pmatrix}.
    \end{equation*}

The trace of $ A $ is zero, and the determinant is $ \det(A) = 0 \cdot 0 - 1 \cdot (-P'(w)) = P'(w) $.

    And the $ (a_{1}, a_{2}) $-subsystem has the coefficient matrix

\begin{equation*}
    B=
    \begin{pmatrix}
        -(a_{1}+a_{2}+w) & -(a_{1}+2w) \\
        a_{2}+2w & a_{1}+a_{2}+w
    \end{pmatrix}.
\end{equation*}

Again, the trace of $ B $ is zero, with determinant $ \det(B) = 3w^{2} - (a_{1}^{2} + a_{1}a_{2} + a_{2}^{2}) $. which is equivalent to $ P'(w) $. 

The principal symbol of the full system is therefore

\begin{equation}
    \sigma(\xi) = 
    \begin{pmatrix}
        \xi_{x}\mathbb{1}_{2} - \xi_{y}A & 0 \\
        0 & \xi_{x}\mathbb{1}_{2} - \xi_{y}B
    \end{pmatrix}.
\end{equation}

Thus, 

\begin{equation}
    \det(\sigma(\xi)) = (\xi_{x}^{2} + P'(w)\xi_{y}^{2})^{2}.
\end{equation}

This is non-zero for every $ \xi \in \mathbb{R}^{2} \setminus \{0\} $ if and only if $ P'(w) > 0 $. Hence, the system is elliptic precisely on this region, and ellipticity degenerates as $ P'(w) \to 0 $.

\end{proof}

\begin{remark}

    We add a remark on this condition for ellipticity. Recall that 

    \begin{equation*}
        P'(w) = (w+a_{1})(w+a_{2}) + (w+a_{1})(w+a_{3}) + (w+a_{2})(w+a_{3}).
    \end{equation*}

    On the admissible domain, $ w \geq -\min\{a_{1}, a_{2}, a_{3}\} $, and so we have $ w + a_{j} \geq 0 $ for every $ j $. Therefore, $ P'(w) \geq 0 $ automatically. Moreover, $ P'(w) = 0 $ precisely when two of the quantities $ w+a_{j} $ vanish. This means that the ellipticity degenerates precisely at points where at least two of the three radii $ |z_{j}| $ vanish.
    
\end{remark}

\subsection{An Almost- Conservation Law}

A useful question is whether the simplest $ S_{3} $-invariant combinations of the $ a_{j} $ functions satisfy conservation laws, with the purpose of seeking a potential formulation. Whilst we shall see that, in this specific search, a conservation law is not obvious, it is illuminating to see that it `almost exists', with a clear obstruction term. Further, this obstruction term is what vanishes when later consider special cases of this construct.

\begin{definition}[$ S_{3} $ And anti-$ S_{3} $-invariant combinations ]\label{DEF:S3AndAntiS3CombinationsofaSector}

    Define
\begin{align*}
    A = a_{1}a_{2}+a_{1}a_{3}+a_{3}a_{2}, \\
    B = a_{1}a_{2}a_{3}, \\
    C = (a_{1}-a_{2})(a_{2}-a_{3})(a_{3}-a_{1}).
\end{align*}

Under the natural action of $ S_{3} $ by permutation of the indices, $ A $ and $ B $ are invariant, while $ C $ is alternating. 
    
\end{definition}

These combinations are particularly natural, as the invariant polynomial $ P(w) $, and the derivative of this polynomial with respect to $ w $, $ P'(w) $ can be written in terms of $ A, B, C $, and what's more, they satisfy a discriminant identity:

\begin{align*}
    P(w) &= w^{3}+Aw+B, \\
    P'(w) &= 3w^{2}+A, \\
    C^{2} &= -4A^{3}-27B^{2}. 
\end{align*}

We now establish two identities satisfied by these combinations, which take the form of conservation laws, plus additional obstruction terms.

\begin{lemma}[Almost-conservation laws]\label{LEMMA:AlmostConservationLaws}

Let $ U \subseteq \mathbb{R}^{2} $ be open, and suppose that $ a_{1}, a_{2}, u, v, w $ are solutions of \eqref{EQ:MainPDESystem} on $ U $, with $ a_{3} = -a_{1}-a_{2} $. 

Then the functions $ A, B, C $ defined in Definition~\ref{DEF:S3AndAntiS3CombinationsofaSector} satisfy

\begin{align}
    \label{EQ:ConservationLawCayleyOne}
    &(A + \frac{3}{2}w^{2})_{x} = -\frac{1}{3}C_{y} - 6wvu_{y}, \\
    \label{EQ:ConservationLawCayleyTwo}
    &C_{x} = 3AA_{y} - 9wB_{y}.
\end{align}

\end{lemma}

\begin{proof} 

The proof is a matter of computation.

Differentiating $ A $ with respect to $ x $, and substituting in the PDEs for $ a_{1, x} $ and $ a_{2, x} $ from \eqref{EQ:MainPDESystem}, we obtain that

\begin{align*}
    A_{x} = &(2a_{1}^{2}+2a_{1}a_{2}-a_{2}^{2}-3a_{2}w)a_{1, y} + \\
    &(a_{1}^{2}-2a_{1}a_{2}-2a_{2}^{2}+3a_{1}w)a_{2, y}.
\end{align*}

Similarly, differentiating $ C $ with respect to $ y $ gives

\begin{align*}
    \frac{1}{3}C_{y} =
&(-2a_{1}^{2}-2a_{1}a_{2}+a_{2}^{2})a_{1, y} - \\
&(a_{1}^{2}-2a_{1}a_{2}-2a_{2}^{2})a_{2, y}.
\end{align*}

Combining these expressions yields

\begin{equation*}
     A_{x} = - \tfrac{1}{3}C_{y} + 3w(-a_{1, y}a_{2} + a_{1}a_{2, y}).
\end{equation*}

Using the defining relation for $ w $, the right-hand-side becomes

\begin{equation*}
    3w(a_{1}a_{2,y}-a_{2}a_{1,y})
=-3ww_{x}+\frac{6wvv_{x}}{P'(w)} = -\left(\frac{3}{2}w^{2}\right)_{x}-6wvu_{y}.
\end{equation*}

This proves the first identity. 

For the second identity, substituting the equations for $ a_{1, x} $ and $ a_{2, x} $ into the derivative of $ C $ gives that

\begin{align*}
C_{x}
=&(2a_{1} + a_{2})
\left(-3A+9wa_{2}\right)a_{1,y}\\
&+(a_{1}+2a_{2})
\left(-3A+9wa_{1}\right)a_{2,y}.
\end{align*}

Noting that

\begin{align*}
A_{y}
&=-(2a_{1}+a_{2})a_{1,y}
  -(a_{1}+2a_{2})a_{2,y},\\
B_{y}
&=-a_{2}(2a_{1}+a_{2})a_{1,y}
  -a_{1}(a_{1}+2a_{2})a_{2,y},
\end{align*}

we conclude that

\begin{equation*}
    C_{x} = 3AA_{y}-9wB_{y},
\end{equation*}

as required.

\end{proof}

Thus, the immediate symmetric candidate does not seem to yield an exact conservation law in this particular case; its failure is measured explicitly by the obstruction term. There are, however, special cases:  in \eqref{EQ:ConservationLawCayleyOne}, the obstruction term when written as $ 3w(-a_{1, y}a_{2}+a_{1}a_{2, y}) $ vanishes when the $ a_{j} $ are constant, which we shall understand in \S~\ref{SECT:SLs} as corresponding to $ \mathrm{U}(1)^{2} $-invariant special Lagrangian $ 4 $-folds. The obstruction term, when written in the form $ 6wvu_{y} $ vanishes when $ v = v_{0} $ and $ u = u_{0} $  are constant, which we shall understand in \S~\ref{SECT:CoassociativeConstrution} as corresponding to $ \mathrm{U}(1)^{2} $-invariant coassociative $ 4 $-folds in $ \mathbb{R}^{7} $, having identified $ \mathbb{R}^{8} = \mathbb{R} \oplus \mathbb{R}^{7} $.

\subsection{Cayley $ 4 $-Folds Along the Degenerate Orbit}

We conclude the section on Cayley $ 4 $-folds by considering the geometry of $ \mathrm{U}(1)^{2} $-invariant Cayley $ 4 $-folds near points at which the $ \mathrm{U}(1) $-orbit degenerates to a $ \mathrm{U}(1) $-orbit. In particular, we are interested in whether a prescribed real-analytic locus of such orbits can be contained in two distinct non-singular $ \mathrm{U}(1)^{2} $-invariant Cayley $ 4 $-folds.

We consider the locus where $ z_{1} = z_{2} = 0 $, and $ z_{3} \neq 0 $. At these points, the $ \mathrm{U}(1)^{2} $-orbit one-dimensional. In terms of our coordinate system, this corresponds to $ w+a_{1} = w+a_{2} = 0 $ and $ w + a_{3} = |z_{3}|^{2} > 0 $, and necessarily, $ v = y = 0 $.

Let $ I \subseteq \mathbb{R} $ be an open interval, and let $ f: I \to (0, \infty) $ and $ \bar{u}: I \to \mathbb{R} $ be real-analytic functions. We prescribe the initial data $ |z_{3}|^{2} = f(x) $, and $ u(x, 0) = \bar{u}(x) $. These conditions imply that $ w = \frac{1}{3}f(x) $, $ A =  -\frac{1}{3}f^{2}(x) $, and $ C = 0 $.

The prescribed locus in $ \mathbb{C}^{4} $ is the real-analytic surface

\begin{equation*}
    \gamma = \bigl\{ \bigl(0, 0, \sqrt{f(x)}e^{i \theta}, x+i\bar{u}(x) \bigr): x\in I,\ \theta \in \mathbb{R}/2\pi\mathbb{Z}\bigr\}.
\end{equation*}

Thus, $ \gamma $ is the union of $ \mathrm{U}(1) $-orbits over a real-analytic curve in the orbit space.

We further expect the local solutions to exhibit a $ \mathbb{Z}_{2} $-equivariance under reflection across the degeneracy $ y = 0 $, with the proposed involution

\begin{equation*}
    (x, y, u, v, a_{1}, a_{2}, a_{3}, w) 
\longmapsto
(x, -y, u, -v, a_{2}, a_{1}, a_{3}, w).
\end{equation*}

This corresponds to $ A $ and $ B $ being even, while $ C $ is odd.

The question is then whether there exist two distinct non-singular $ \mathrm{U}(1)^{2} $-invariant Cayley $ 4 $-folds $ N_{+} $ and $ N_{-} $ containing $ \gamma $, whose reduced descriptions correspond locally to solutions of \eqref{EQ:MainPDESystem} on the regions $ y > 0 $ and $ y < 0 $, respectively. In particular, we would like to understand whether the two branches can intersect precisely along $ \gamma $, giving rise to a possibly non-singular union $ N = N_{+} \cup N_{-} $.

The principal analytic difficulty arises from the degeneration of the reduced coordinate system along $ y = 0 $. In particular, the functions $ u, v, a_{1}, a_{2}, w $, and hence the combinations $ A, B, C $, are not in general smooth across $ y = 0 $. Such a loss of regularity in the reduced variables does not, however, necessarily preclude smoothness of the corresponding submanifolds in $ \mathbb{C}^{4} $.

In \cite{joyce2005u3}, Joyce is able to overcome related difficulties in smoothness due to a proof of the existence of singular solutions. A complementary approach is provided by Bryant \cite{bryant2006so}, who establishes existence results for non-singular special Lagrangian submanifolds with prescribed real-analytic fixed loci.

It would be of interest to determine whether an analogous theory can be developed for our Cayley construct, overcoming the current issue of the loss of smoothness across $ y = 0 $. One possible approach is to seek coordinates adapted to the non-principal orbit structure, in which the equations governing smooth extensions across $ \gamma $ become more tractable. We leave this for future investigation.

\section{$ \mathrm{U(1)}^{2} $-Invariant Coassociative $ 4 $-Folds in $ \mathbb{R}^{7} $}\label{SECT:CoassociativeConstrution}

In this section, we consider the first special case of Theorem~\ref{THM:MainTheoremRegular}, which is a construction for $ \mathrm{U(1)}^{2} $-invariant coassociative $ 4 $-folds which live in the affine plane $ x_{8} = u_{0} $, on $ \mathbb{R}^{7} \subset \mathbb{R}^{8} $. Taking that $ u = u_{0} \in \mathbb{R} $, first two equations in \eqref{EQ:MainPDESystem} imply that $ v_{y} = 0 $ and $ v_{x} = 0 $. Hence, $ v = v_{0} $ is locally constant, which we can view as like a moment map equation for $ \phi_{0} $.

\begin{theorem}[Construction of a class of $ \mathrm{U(1)}^{2} $-invariant coassociative submanifolds]\label{THM:ConstructionOfCoassociativeSubmanifolds}

    Let $ (x, y) $ be coordinates on a smooth domain $ \Sigma \subset \mathbb{R}^{2} $, and let $ a_{1}, a_{2}, a_{3}: \Sigma \to \mathbb{R} $ be continuous functions satisfying $ a_{1}+a_{2}+a_{3} = 0 $. Fix constants $ u_{0} = \mathrm{Im}(z_{4}) $ and $ v_{0} = \mathrm{Re}(iz_{1}z_{2}z_{3}) $. 
    
    Define the degeneracy set

    \begin{equation*}
        S := \{ (x, y) \in \Sigma: y = v_{0} = 0, \min(a_{j}(x, y)) \text{ has multiplicity} \geq 2 \},
    \end{equation*}

    and set $ \Sigma^{*} := \Sigma \backslash S $. For $ (x, y) \in \Sigma^{*} $, let $ w = w(x, y) $ denote the admissible root of $ (w+a_{1})(w+a_{2})(w+a_{3}) = v_{0}^{2} + y^{2} $.

    Define the $ \mathrm{U(1)}^{2} $-invariant submanifold $ N \subset \mathbb{C}^{4} $ by

    \begin{equation*}
    \begin{aligned}
        N := \{ (z_{1}, z_{2}, z_{3}, z_{4}) \in H_{u_{0}} \, | \, &|z_{1}|^{2} - |z_{3}|^{2} = a_{1}(x, y) - a_{3}(x, y), \\
        &|z_{2}|^{2}-|z_{3}|^{2} = a_{2}(x, y) - a_{3}(x, y), \\
        &iz_{1}z_{2}z_{3} = v_{0} + iy, \\
        &z_{4} = x+iu_{0}, \quad (x, y) \in \Sigma^{*}\}.
    \end{aligned}
    \end{equation*}

    If the PDE system 

    \begin{equation}\label{EQ:CoassociativePDESystem}
   \begin{cases}
    a_{1,x} = -\bigl((a_{1}+a_{2}+w)a_{1,y}+(a_{1}+2w)a_{2,y}\bigr), \\
a_{2,x} = (a_{2}+2w)a_{1,y} + (a_{1}+a_{2}+w)a_{2,y},
\end{cases}
\end{equation}
    
    holds on $ \Sigma^{*} $, then $ N $ is a smooth $ \mathrm{U(1)}^{2} $-invariant coassocaitive $ 4 $-fold on the free locus.
    
\end{theorem}

The proof of this theorem is identical to that of Theorem~\ref{THM:MainTheoremRegular}, except we now take $ u = u_{0} $ and $ v = v_{0} $ constant.

\subsection{Analysis of $ \mathrm{U(1)}^{2} $-Invariant Coassociative $ 4 $-Folds}

The coassociative submanifold sector tells a different analytical story than the most general Cayley case, in particular, a potential function now exists. Indeed, the behaviour of these submanifolds is comparable to the $ \mathrm{U(1)} $-invariant special Lagrangian $ 3 $-folds of \cite{joyce2005u1}, letting $ a \leftrightarrow v_{0} $. 

\begin{theorem}[Potential formulation for coassociative submanifolds]

    Let $ \Sigma^{*} \subset \mathbb{R}^{2} $ be as in Theorem~\ref{THM:ConstructionOfCoassociativeSubmanifolds}, let $ (a_{1}, a_{2}, a_{3}, w) $ be data satisfying \eqref{EQ:CoassociativePDESystem}, and let $ U \subset \Sigma^{*} $ be open and simply-connected. Then on  $ U $,

\begin{equation}\label{EQ:CoassociativeManifoldsConservationLaw}
    (A + \frac{3}{2}w^{2})_{x} = -\frac{1}{3}C_{y}.
\end{equation}

The one-form

\begin{equation}
    \alpha := -\frac{1}{3} C dx + \bigl( A + \frac{3}{2}w^{2} \bigr)dy
\end{equation}

is closed on $ U $, and is exact. Therefore, there exists a potential function $ g: U \to \mathbb{R} $, unique up to an additive constant, with

\begin{equation*}
    g_{y} = (A+\frac{3}{2}w^{2}), \quad  g_{x} = -\frac{1}{3}C.
\end{equation*}

The function $ w $ is a root of the equation

\begin{equation}\label{EQ:gxgyEQForCoassociative}
    9g_{x}^{2} = -4(g_{y}-\frac{3}{2}w^{2})^{3} - 27(v_{0}^{2}+y^{2}-wg_{y} +\frac{1}{2}w^{3})^{2},
\end{equation}

On an open subset of $ U $, where \eqref{EQ:gxgyEQForCoassociative} admits this admissible branch $ w = W(g_{x}, g_{y}, y) $ smoothly, the PDE system \eqref{EQ:CoassociativePDESystem} is equivalent to the scalar quasilinear equation

\begin{equation}\label{EQ:SecondOrderPotentialEQ}
    g_{xx} + P'(w)g_{yy} - 6wy = 0.
\end{equation}
    
\end{theorem}

\begin{proof}

By Lemma~\ref{LEMMA:AlmostConservationLaws}, 

\begin{equation*}
    (a_{1}a_{2}+a_{2}a_{3} + a_{3}a_{1}+\frac{3}{2}w^{2})_{x} = -\frac{1}{3}((a_{1}-a_{2})(a_{2}-a_{3})(a_{3}-a_{1}))_{y} + \frac{6wvv_{x}}{P'(w)}.
\end{equation*}

Since we are working with $ u = u_{0} $, $ v = v_{0} $ are constant, $ v_{x} = 0 $ and \eqref{EQ:CoassociativeManifoldsConservationLaw} follows.

Define $ \alpha $ on $ U $ as above. Then $ d\alpha = 0 $ by \eqref{EQ:CoassociativeManifoldsConservationLaw}, so since $ U $ is simply-connected, $ \alpha = dg $ for some $ g \in C^{2}(U) $, unique up to an additive constant, giving $ g_{x}, g_{y} $ as claimed.

Equating the coefficients of $ dx $ and $ dy $ in $ \alpha = dg $ then gives us that $ g_{x} = -\tfrac{1}{3}((a_{1}-a_{2})(a_{2}-a_{3})(a_{3}-a_{1})), g_{y} = (a_{1}a_{2}+a_{2}a_{3}+a_{3}a_{1}+\tfrac{3}{2}w^{2}) $. Thus, our claimed equation follows by substituting this into \eqref{EQ:MainPDESystem}.
    
\end{proof}

\begin{remark}

    The nature of the potential equation for the coassociative equation suggests the natural question: \textit{can one prove an existence and uniqueness theorem for $ \mathrm{U(1)}^{2} $-invariant coassociative $ 4 $-folds in $ \mathbb{R}^{7} $?} This lends itself to future studies.
    
\end{remark}

\section{Construction and Analysis of $ \mathrm{U(1)}^{2} $-Invariant Special Lagrangian $ 4 $-Folds in $ \mathbb{R}^{8} $}\label{SECT:SLs}

This section treats the final specialisation of Theorem~\ref{THM:MainTheoremRegular}, namely the construction of $ \mathrm{U(1)}^{2} $-invariant special Lagrangian four-folds in $ \mathbb{C}^{4} $. In this sector, the functions $ a_{1}, a_{2}, a_{3} $ are constant, while $ u $ and $ v $ remain functions of the reduced coordinates. The setting and theorem is then identical to the $ n = 4 $ case in \cite{MyPaper}, and serves as a generalisation of \cite{joyce_special_2002}.

\begin{theorem}[A construction of $ \mathrm{U(1)}^{2} $-invariant special Lagrangian $ 4 $-folds]

Let $ (x, y) \in \Sigma \subset \mathbb{R}^{2} $ be coordinates in a smooth domain, and fix parameters $ (a_{1}, a_{2}, a_{3}) \in \mathbb{R}^{3} $, satisfying $ a_{1}+a_{2}+a_{3} = 0 $. For each $ (x, y) \in \Sigma $, let $ w = w(x, y) $ denote the unique admissible root of $ (w+a_{1})(w+a_{2})(w+a_{3}) = v^{2}+y^{2} $, with $ w \geq - \min_{j}(a_{j}) $.

 Define the degeneracy set

 \begin{equation*}
        S := \{ (x, y) \in \Sigma: y = v(x, y) = 0, \min(a_{j}(x, y)) \text{ has multiplicity} \geq 2 \},
    \end{equation*}

 and $ \Sigma^{*} = \Sigma \backslash S $ to be the smooth domain minus the degenerate points, where $ u, v \in C^{0}(S) $ on $ \Sigma^{*} $. Let $ w = w(x, y) $ be the unique admissible root of $ (w+a_{1})(w+a_{2})(w+a_{3}) = v^{2}+y^{2} $ on $ \Sigma^{*} $. 

  Define

    \begin{equation*}
    \begin{aligned}
        N := \{ z \in \mathbb{C}^{4} \, | \, |z_{1}|^{2} - |z_{3}|^{2} &= a_{1} - a_{3}, \\
        |z_{2}|^{2} - |z_{3}|^{2} &= a_{2} - a_{3}, \\
        iz_{1}z_{2}z_{3} &= v(x, y) + iy, \\
        z_{4} &= x+iu(x, y), \quad (x, y) \in \Sigma^{*} \}.
    \end{aligned}
    \end{equation*}

    If the PDE system 

    \begin{equation}\label{EQ:SLPDESystem}
        \begin{cases}
            u_{x} = v_{y}, \\
            v_{x} = -P'(w)u_{y},
        \end{cases}
    \end{equation}

    holds on $ \Sigma^{*} $, then $ N $ is a smooth $ \mathrm{U(1)}^{2} $-invariant special Lagrangian $ 4 $-fold on the free locus.
    
\end{theorem}

The proof of this theorem is again identical to the proof of Theorem~\ref{THM:MainTheoremRegular}, except that we now take $ a_{1}, a_{2}, a_{3} \in \mathbb{R} $ are constant.

If $ (u, v) $ satisfy \eqref{EQ:SLPDESystem}, then the condition $ u_{x} = v_{y} $ implies the local existence of a potential function $ f $ such that $ f_{y} = u $ and $ f_{x} = v $, which leads us to the following proposition.

\begin{proposition}[Potential formulation for the special Lagrangian PDE system]

Let $ S \subset \mathbb{R}^{2} $ be a simply-connected domain, and let $ u, v \in C^{1}(S) $ satisfy \eqref{EQ:SLPDESystem}. Then there exists $ f \in C^{2}(S) $, unique up to the addition of a constant, such that

\begin{equation*}
    f_{x} = v, \quad f_{y} = u.
\end{equation*}

Moreover, $ f $ satisfies the second-order equation

\begin{equation}\label{EQ:SecondOrderSLEQ}
    f_{xx} + P'(w)f_{yy} = 0,
\end{equation}

where $ w = w(f_{x}, y) $ is the unique real solution of $ P(w) = f_{x}^{2} + y^{2} $, satisfying $ w \geq - \min(a_{j}) $. 

Conversely, any function $ f \in C^{2}(S) $ solving \eqref{EQ:SecondOrderSLEQ} yields a solution $ (u, v) = (f_{y}, f_{x}) $ of \eqref{EQ:SLPDESystem}.
    
\end{proposition}

Equation \eqref{EQ:SecondOrderSLEQ} is a second-order quasilinear elliptic equation wherever $ P'(w) > 0 $, and becomes degenerate precisely along the locus $ P'(w) = 0 $. 

A generalisation of the existence and uniqueness theory developed in \cite{joyce2005u1} follows, as outlined in \cite{MyPaper}, in the non-singular case. Extending this theory to the singular case is a separate problem. It would require an analogue of Joyce's weak-solution and singularity analysis, which we do not pursue here.

\subsection{Example of a $ \mathrm{U(1)}^{2} $-Invariant Special Lagrangians}

We conclude this section by considering the key example that Harvey and Lawson \cite{harvey_calibrated_1982} provided.

\begin{theorem}

    Let $ M_{c} $ denote the locus in $ \mathbb{C}^{4} $ of the equations

    \begin{equation}\label{EQ:HarveyLawsonMomentMaps}
        |z_{1}|^{2} - |z_{2}|^{2} = c_{2}, \quad |z_{1}|^{2} - |z_{3}|^{2} = c_{3}, \quad |z_{1}|^{2} - |z_{4}|^{2} = c_{4},
    \end{equation}

    and

    \begin{equation}\label{EQ:HarveyLawsonRealEquation}
        \mathrm{Re}(z_{1}z_{2}z_{3}z_{4}) = c.
    \end{equation}

    The smooth locus of $ M_{c} $ is special Lagrangian. When the defining equations have full rank, $ M_{c} $ is a smooth special Lagrangian submanifold.
     
\end{theorem}

Let's consider this in our setting. Recall that $ z_{4} = x+iu $ and $ iz_{1}z_{2}z_{3} = v+iy $ are our invariants. Set that $ \kappa_{1} = c_{4} $, $ \kappa_{2} = c_{4} - c_{2} $, and $ \kappa_{3} = c_{4} - c_{3} $, and let

\begin{equation*}
    \tilde{\kappa}=\frac{1}{3}
(\kappa_{1}+\kappa_{2}+\kappa_{3}).
\end{equation*}

The moment-map constraints imply

\begin{equation*}
    |z_{j}|^{2}=x^{2}+u^{2}+\kappa_{j},\qquad j=1,2,3.
\end{equation*}

Thus, in our earlier notation, we have

\begin{equation*}
    w=x^{2}+u^{2}+\tilde{\kappa},\qquad
a_{j}=\kappa_{j}-\tilde{\kappa}.
\end{equation*}

In particular, the parameters $ a_{j} $ are constant, and satisfy the normalisation constraint $ a_{1} + a_{2} + a_{3} = 0 $.

Since

\begin{equation*}
    z_{1}z_{2}z_{3}z_{4}
=(y-iv)(x+iu)
=xy+uv+i(uy-xv),
\end{equation*}

the defining equations reduce to

\begin{equation}
\label{eq:HL-Reduced}
\begin{cases}
v^{2}+y^{2}=\displaystyle\prod_{j=1}^{3}
(x^{2}+u^{2}+\kappa_{j}),\\[4pt]
xy+uv=c.
\end{cases}
\end{equation}

Let $ P(w)=\prod_{j=1}^{3}(w+a_{j}) $, as before. The Jacobian determinant of the two equations in \eqref{eq:HL-Reduced}, with respect to $ (u, v) $, is

\begin{equation*}
    -2\bigl(v^{2}+u^{2}P'(w)\bigr).
\end{equation*}

Thus, wherever $ D:= v^{2} + u^{2} + 2P'(w) > 0 $, the Implicit Function Theorem allows us to express $ u $ and $ v $ locally as smooth functions of $ (x, y) $.

Differentiating \eqref{eq:HL-Reduced} gives

\begin{equation*}
    u_{x}=v_{y}=-\frac{xuP'(w)+yv}{D},
\end{equation*}

and

\begin{equation*}
    u_{y}=\frac{uy-vx}{D},\qquad
v_{x}=-P'(w)\frac{uy-vx}{D}.
\end{equation*}

Hence,

\begin{equation*}
    u_{x}=v_{y},\qquad v_{x}=-P'(w)u_{y},
\end{equation*}

which is exactly our special Lagrangian system.

We can also pursue an algebraic description of the solutions. On the locus $ u \neq 0 $, we have $ v=\tfrac{c-xy}{u} $. Substituting into \eqref{eq:HL-Reduced}, we obtain 

\begin{equation}
\label{eq:HL-quartic}
u^{2}\left[
\prod_{j=1}^{3}(x^{2}+u^{2}+\kappa_{j})-y^{2}
\right]-(c-xy)^{2}=0.
\end{equation}

This is a quartic polynomial in $ u^{2} $, and each positive root which satisfies $ x^{2} + u^{2} + \kappa_{j} \geq 0 $ for all $ j $ determines the candidate branches. Where $ D > 0 $, the branches are locally smooth.

This recovers the Harvey-Lawson family without our $ \mathrm{U}(1)^{2} $-invariant special Lagrangian construction.

\section*{Conclusion}

In this paper, we derived a non-linear, Cauchy-Riemann system of partial differential equations characterising $ \mathrm{U}(1)^{2} $-invariant Cayley $ 4 $-folds in $ \mathbb{C}^{4} $, with ellipticity away from the degeneracy locus. We further considered two distinguished sectors: the coassociative $ 4 $-folds, where we take $ u = u_{0} $ and $ v = v_{0} $ as constant, and the special Lagrangian $ 4 $-folds, where we take $ a_{1}, a_{2}, a_{3} $ as constant.

In the most general Cayley setting, we have considered the $ S_{3} $ and anti-$ S_{3} $-invariant combinations of the moment-map parameters, and derived 'almost conservation laws', by which we mean a conservation law plus an obstruction term. By contrast, the special Lagrangian and coassociative reductions admit potential formulations.

There are several natural questions which arise from these results:

\begin{itemize}

    \item \textit{With the potential formulation for the class of $ \mathrm{U}(1)^{2} $-invariant coassociative $ 4 $-folds in $ \mathbb{R}^{7} $, can we prove the existence and uniqueness of solutions?}

    \item \textit{Despite not finding a potential formulation from the simplest $ S_{3} $-invariant combinations of the $ a_{j} $-functions, could we expand our search to seek a potential from more complicated combinations, or by using alternate coordinate systems? 
    If this is not possible, can we prove existence and uniqueness for the family of $ \mathrm{U}(1)^{2} $-invariant Cayley $ 4 $-folds in $ \mathbb{R}^{8} $ by setting up the problem as a boundary elliptic problem, potentially using Atiyah-Patodi-Singer boundary conditions \cite{AtiyahPatodiSingerBoundaryConditionsI, AtiyahPatodiSingerBoundaryConditionsII, AtiyahPatodiSingerBoundaryConditionsIII}, specifying half of the $ u, v, a_{1}, a_{2} $ data on the boundary?}

    \item \textit{Motivated by \cite{LotayCoassociativeCones, lotay2007calibrated},  we can ask whether the reduced equations admit $ \mathrm{U}(1)^{2} $-invariant coassociative or Cayley cones whose links are diffeomorphic to $ T^{3} \subset S^{6} $ or $ S^{7} $? In particular, can suitable homogeneity conditions reduce the system to equations admitting periodic systems?}
    
\end{itemize}

These questions suggest several directions for further study of the analytic and geoemtric properties of the reduced Cayley equations.

\renewcommand{\bibsection}{%
  \section*{References}%
}

\bibliography{references}
\bibliographystyle{plainnat}

\end{document}